\documentclass{article}
\usepackage[utf8]{inputenc}

\usepackage[legalpaper, margin=1in]{geometry}
\usepackage{amsmath,amssymb,times,bm}
\usepackage{stmaryrd}
\usepackage{enumerate}
\usepackage{relsize}
\usepackage{graphicx}
\usepackage{subfigure}
\usepackage{mathrsfs}
\usepackage{latexsym}
\usepackage{amsfonts}
\usepackage{pifont}
\usepackage{nicematrix}
\usepackage{pb-diagram}
\usepackage{setspace}
\usepackage{url} 
\usepackage{hyperref}
\usepackage{algorithm}
\usepackage{algorithmic}
\usepackage{color}
\usepackage{multirow}
\usepackage{mathtools}
\usepackage{appendix}
\usepackage{svg}
\usepackage{lineno}
\usepackage{tikz}
\usepackage{amsmath}
\usepackage{amsfonts}
\usepackage{amssymb}
\usepackage{amsthm}
\usepackage{caption}
\usepackage{subcaption}
\usepackage{arydshln}
\usepackage{verbatim}
\usepackage{bbold}
\usepackage[normalem]{ulem}
\allowdisplaybreaks
\NewDocumentCommand{\tens}{t_}
 {%
  \IfBooleanTF{#1}
   {\tensop}
   {\otimes}%
 }
\usepackage{pgfplots}
\usetikzlibrary{shapes,arrows, arrows.meta,calc,tikzmark}
\pgfplotsset{compat=newest}

\pgfplotsset{compat=1.16}
\usetikzlibrary{shapes,positioning,intersections,quotes}
\tikzset{
    cross/.pic = {
    \draw[rotate = 45] (-#1,0) -- (#1,0);
    \draw[rotate = 45] (0,-#1) -- (0, #1);
    }
}

\newcounter{anote} 

\newtheorem{weakform}{Weak form}

\title{Two Nitsche-based mixed finite element discretizations for the seepage problem in Richards' equation}

\author{
Federico Gatti$^{(a),*}$  
Andrea Bressan$^{(a)}$ 
Alessio Fumagalli$^{(d)}$
Domenico Gallipoli$^{(b)}$ \\
Leonardo Maria Lalicata$^{(b)}$ 
Simone Pittaluga$^{(c)}$ 
Lorenzo Tamellini$^{(a)}$ 
}

\date{}

\def\Gf{\Gamma_{\text{flux}}}
\def\Gb{\Gamma_{\text{bot}}}
\def\Gt{\Gamma_{\text{top}}}
\def\Gl{\Gamma_{\text{lat}}}
\def\Gbh{\Gamma_{\text{bot},h}}
\def\Gth{\Gamma_{\text{top},h}}
\def\Glh{\Gamma_{\text{lat},h}}
\def\Gfh{\Gamma_{\text{flux},h}}
\def\GthN{\Gamma_{\text{top},h}^{\text{no-hyb}}}
\def\GthH{\Gamma_{\text{top},h}^{\text{hyb}}}

\def\vq{\mathbf{q}}
\def\vn{\mathbf{n}}
\def\vN{\mathbf{N}}

\def\vp{\mathbf{p}}
\def\vv{\mathbf{v}}
\def\vA{\mathbf{A}}
\def\vB{\mathbf{B}}
\def\vC{\mathbf{C}}
\def\vD{\mathbf{D}}
\def\vE{\mathbf{E}}
\def\vF{\mathbf{H}^{(\vq)}}
\def\vG{\mathbf{H}^{(\psi)}}
\def\vH{\mathbf{H}}
\def\vI{\mathbf{I}}
\def\vV{\mathbf{V}}

\def\tK{\mathbf{K}}
\def\RT{\mathbb{RT}_0}
\def\P{\mathbb{P}_0}

\begin{document}

\maketitle

\begin{center}
{ 
    \small $^{a}$ Consiglio Nazionale delle Ricerche -- Istituto di Matematica Applicata e Tecnologie Informatiche ``E. Magenes''  (CNR-IMATI), 27100, Pavia, Italy\\
    {
        \tt
        federico.gatti@imati.cnr.it
    }\\
    \small $^{b}$Department of Civil, Chemical and Environmental Engineering, University of Genoa, 16145, Genoa, Italy\\
    \small $^{c}$Consiglio Nazionale delle Ricerche -- Istituto di Matematica Applicata e Tecnologie Informatiche ``E. Magenes'' (CNR-IMATI), 16145, Genoa, Italy\\
    \small $^{d}$MOX, Department of Mathematics, Politecnico di Milano, 20133, Milan, Italy\\
}
\end{center}

\begin{abstract}


This paper proposes two algorithms to impose seepage boundary conditions in the context of Richards' equation
for groundwater flows in unsaturated media. Seepage conditions are non-linear boundary conditions, that can be 
formulated as a set of unilateral constraints on both the pressure head and the water flux at the ground surface,
together with a complementarity condition:
these conditions in practice require switching between Neumann and Dirichlet boundary conditions on unknown portions on the boundary.
Upon realizing the similarities of these conditions with unilateral contact problems in mechanics,
we take inspiration from that literature to propose two approaches:
the first method relies on a strongly consistent penalization term, whereas the second one is obtained by an hybridization approach, in which the value of the pressure on the surface is treated as a separate set of unknowns. 
The flow problem is discretized in mixed form with div-conforming elements so that the water mass is preserved.
Numerical experiments show the validity of the proposed strategy in handling the seepage boundary conditions on geometries with increasing complexity.

\end{abstract}

\noindent
{\bf Keywords}: Richards' equation, \and Seepage problem, \and Nitsche's method, \and Signorini problem, 
\and Mixed finite elements, \and Hybridized finite elements.\\



\section{Introduction}\label{sec1}
One of the difficulties in the simulation of ground
water flows is dealing with changing regimes at ground surface.
Where the soil is not saturated, water infiltrates; conversely where the soil is saturated either water accumulates forming ponds, or it runs off depending on the terrain geometry.
In the following, it is assumed that the ground surface of interest is
a slope so that the formation of ponds is excluded, and the
saturated regions result in the formation of \emph{seepage faces},
i.e. spring-like regions where ground water exits the terrain and becomes surface-water.
The incorrect determination of seepage
faces causes incorrect prediction of key parameters such as pore water
pressures, effective stresses, and shear strengths within soil masses,
and impacts the assessment of the stability of 
slopes~\cite{rulon1985multiple, crosta1999slope, orlandini2015evidence}, embankments~\cite{hirschfeld1973embankment, milligan2003some},
and other geotechnical structures.
Therefore, the correct determination of seepage faces is of paramount importance
for identifying and mitigating potential hazards such as landslides 
and slope failures, which pose substantial risks to infrastructure and lives~\cite{froude2018global}.
The motivation of this work is indeed the detection of rainfall-induced landslide triggering:
to this end, the ground water flow model needs to be coupled with stability analysis algorithms
(either limit equilibrium models such as \cite{bishop1955use, morgenstern1965analysis} or 
stress-based method), and the effects of the landslides can be assessed by landslide runout simulations,
see e.g. \cite{pastor2021depth,he2023mpm,GATTI2024128525,GATTI2024112798}. 

The ground water flow is mathematically modelled by the Richards'
equation, and the different ground surface flow regimes
correspond to different types of boundary conditions: 
the water flux entering the unsaturated region is modeled with a Neumann type boundary condition, 
whereas water runoff on the surface of the saturated region is modeled by imposing that the value of the pressure head be zero, i.e. a homogeneous Dirichlet type boundary condition. 
We remark that we do not know in advance on what portions of the boundary each condition must be imposed,
which will finally imply solving a non-linearity in the final system of equations.


To address the determination of seepage faces in the solution of Richards' equation, several finite element algorithms have been proposed in the geotechnical literature. Starting from the pioneering work of~\cite{freeze1971three}, which presents one of the first three-dimensional finite difference models in the presence of seepage faces, an iterative algorithm for finite element schemes has been proposed in a series of papers, including those by~\cite{rubin1968theoretical, neuman1975finite, cooley1983some}. More recent attempts to generalize the algorithm can be seen in works such as~\cite{scudeler2017examination}.
This algorithm employs an iterative procedure based on the switching between flux-type (Neumann) and head-type (Dirichlet) boundary conditions to determine seepage faces and locate the so-called ``exit points'',
i.e., the locations on a slope (such as the downstream face of an earth dam or embankment) where the phreatic line intersects with the surface, where seepage water exits the earth body.
More in details, the iterative algorithm begins each time step by setting 
all potential seepage faces to a flux-type boundary condition and computing a first-attempt solution. 
Once pressure heads are available, the algorithm: 
a) searches for subregions of the ground surface where the pressure is positive;
b) forces the corresponding pressure degrees of freedom to zero, 
imposing homogeneous Dirichlet boundary conditions;
c) computes a new proposed solution, iterating these three steps until no positive pressure on the surface is detected. 
However, this heuristic procedure has several drawbacks from a mathematical perspective.    
Primarily, it lacks a well-established mathematical theory to address convergence and stability issues, 
and furthermore it may exhibit water mass conservation problems.

In this paper, we propose two new strategies to accurately determine seepage faces, by drawing similarities between the boundary conditions discussed above (to which we will refer in the rest of this work we as \emph{seepage conditions} for brevity) and the Karush-Kuhn-Tucker conditions arising in contact mechanics problems, together with their numerical treatment discussed in the series of works \cite{chouly2013nitsche, chouly2015symmetric, chouly2017overview, chouly2018unbiased}. 
Specifically, we consider a mixed finite element formulation \cite{raviart2006mixed} 
of the Richards' equation, and the first strategy
that we propose consists in imposing the seepage conditions with a ``Nitsche-inspired'' approach, where 
the conditions are imposed in a weak way relying on a penalization term. 
The second strategy that we propose is to move further to a hybridized formulation \cite{brezzifortin1991},
in which the pressure head on the boundary is treated as a separate unknown of the problem.
The cost of adding further degrees of freedom on the potential seepage faces is demonstrated to be advantageous as it frees the computation from the choice of the penalization parameter appearing in the first strategy. 


As already mentioned, the methods that we propose are based on the scheme developed in the series of works \cite{chouly2013nitsche, chouly2015symmetric, chouly2017overview, chouly2018unbiased} 
for the simulation of contact problems in linear elastic equations,
in which a variant of the Nitsche's method is used to impose the contact boundary conditions in a weak way. The Nitsche method has been designed as a valid alternative to the standard Lagrange multiplier approach~\cite{tur2009mortar, brunssen2007contact, papadopoulos1998lagrange} in the field of contact problems and has gained popularity as a fully consistent penalization method.
Later, this contact algorithm has been applied in the framework of fluid-structure interaction problems~\cite{burman2020nitsche, burman2007stabilized, burman2014unfitted, formaggia2021xfem} and in the field of fractured porous media~\cite{berge2020finite}. 
The mixed finite element scheme, together with the control-volume finite elements~\cite{cai1997control} and the finite volume methods~\cite{aavatsmark2002introduction}, is among the classical conservative discretizations used for the approximation of the elliptic equations arising in porous media analysis. 

Through numerical experiments on both simplified and realistic terrain models, we validate the efficacy of our two proposed approaches in imposing the seepage conditions on increasingly complex cases. Furthermore, we conduct a comparative analysis between the proposed discretization approaches to highlight their respective strengths and potentials. In details, our findings are that while both approaches work well,
the hybridized approach is generally to be preferred, since the quality 
of the solution does not depend on the tuning of the penalization parameter.
Moreover, in certain situations the temporal evolution of the
pressure and flux field predicted by the non-hybridized scheme shows
a non-negligible delay, whereas this issue is not affecting 
the hybridized scheme if a suitable relaxation of the seepage conditions
is introduced.

The subsequent sections of this paper are structured as follows: In Section~\ref{sec:mode_equations}, we present the governing equations; Section~\ref{sec:discretization} is dedicated to discussing the discretization methods, encompassing both the mixed finite element discretization for Richards' equation and the our two strategies for the seepage problem. In Section~\ref{sec:examples}, we showcase numerical experiments of increasing geometric complexity, ranging from a simple rectangular soil column to a realistic topography derived from a Digital Terrain Model (DTM). These examples aim to demonstrate the versatility and applicability of our proposed methods across a spectrum of real-world scenarios, beyond academic settings. 
Finally, in Section~\ref{sec:conclusion} we draw some conclusions and future perspectives.

\section{Problem statement}\label{sec:mode_equations}

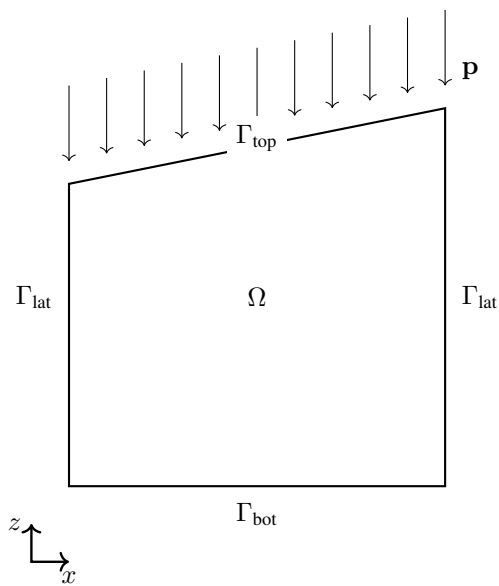
\begin{figure}[ht]
\centering
\begin{tikzpicture}
\draw[black, thick] (0,0)--(5,0)--(5,5)--(0,4)--cycle;
\def\factor{0.15}
\newcommand\y{6}
\foreach \x[evaluate=\x as \y using {.2*\x+5.3}] in {0,0.5,...,5} {
     \draw[->] (\x,\y)--++(0,-1);
}
\node (A) at (-1,-1) {};
\draw[<->,line width=0.9]
    (A)++(0.5,0.5) node[left,scale=1] {$z$} |-++
         (0.5,-0.5) node[below,scale=1] {$x$};

\node[black, right=1pt of {(5,5.5)}, outer sep=2pt,fill=white] {$\vp$};
\node[black, left=1pt of {(0,2.5)}, outer sep=2pt,fill=white] {$\Gl$};
\node[black, right=1pt of {(5,2.5)}, outer sep=2pt,fill=white] {$\Gl$};
\node[black, below=1pt of {(2.5,0)}, outer sep=2pt,fill=white] {$\Gb$};
\node[black, below=1pt of {(2.5,5)}, outer sep=2pt,fill=white] {$\Gt$};
\node[black] at (2.5, 2.5) {$\Omega$};
\end{tikzpicture}
\caption{Sketch of the porous domain $\Omega$. The arrows represent the 
rainfall on the ground surface, that is modeled as a vector field $\vp$.
}
\label{fig:1}
\end{figure}

From the physical point of view the simulation domain $\Omega$ is a 2D vertical section of the terrain, extending from the ground surface to   possibly below the unsaturated zone.
The boundary of the computational domain $\partial \Omega$ is partitioned in three disjoint regions: $\Gt$ is the ground surface (possibly a slope), $\Gb$ is the base of $\Omega$, and $\Gamma^{\text{lat}}$ is the vertical part of the boundary, see Figure \ref{fig:1}. 
We assume that no ponding can occure on $\Gt$ (for example, this is reasonable when $\Gt$ is inclined).
We also assume that there is no water flux through $\Gl$.
More physically realistic conditions depends on the specific application, see e.g. \cite{bianchi2022analysis} for an example in the contest of slope stability.
The main topic of this article is modelling the water flux on $\Gt$, where the flux cannot exceed a given precipitation (rainfall) and, at the same time, pressure cannot be positive 
(the seepage conditions already discussed in the introduction).

\subsection{Richards' equation}\label{sec:rich_1}
The unsaturated flow in $\Omega$ for the time interval $[0,t_{\text{fin}}]$  is modelled by Richards' equation that describe the relation between two fields: the Darcy velocity $\vq(\mathbf{x},t) : \Omega\times[0,t_{\text{fin}}] \rightarrow \mathbb{R}^d$ and the pressure head $\psi(\mathbf{x},t) : \Omega\times[0,t_{\text{fin}}] \rightarrow
\mathbb{R}$, specifically:
\begin{equation}\label{eq:model1:1}
  \begin{aligned}
\tK^{-1}(\psi)\vq + \nabla\psi + \nabla z&= \mathbf{0} &\text{in } &\Omega\times(0,t_{\text{fin}}), \vspace{.2cm}\\
\partial_t\theta(\psi) + \nabla\cdot\vq &= 0 &\text{in } &\Omega\times(0,t_{\text{fin}}),\vspace{.2cm} \\
\psi(\mathbf{x},0) &= \psi_0 &\text{in } &\Omega. 
  \end{aligned}
\end{equation}
The quantity $\nabla z$ is the unit vector $\mathbf{e}_z$ and it identifies the direction opposite to gravity. The permeability tensor $\tK(\psi)$ and the water retention curve $\theta = \theta(\psi)$ are given by Van~Genuchten constitutive laws, see~\cite{van1980closed}. 
In detail, letting $[\cdot]_{-}= \min(0,\cdot)$,    
defining
\[
\hat\theta(\psi):= \left[ 1+ (-\alpha[\psi]_-)^m \right]^{-\frac{m-1}{m}} 
\]
and noting that $0< \hat\theta(\psi) \leq 1$, the water content is taken as
\[
\theta(\psi)=\theta_R+\hat\theta(\psi)(\theta_S-\theta_R),
\]
where $\theta_R$ and $\theta_S$ are respectively the residual and saturated water-content, and the permeability is assumed to be isotropic, i.e., $\tK = K(\psi)\vI$, where $\vI$ is the identity matrix, 
\[ 
K(\psi)=K_S\, \hat\theta(\psi)^{1/2}\Big(1-(1-\hat\theta(\psi)^\frac{m}{m-1})^{\frac{m-1}{m}} \Big)^2
\]
where $K_S$ is the saturated permeability.

\subsection{Boundary conditions for the seepage problem}\label{sec:boundary_cond}

The system \eqref{eq:model1:1} is closed given proper boundary data.
To fix the ideas, we consider that there is no flux through $\Gl$ and $\Gb$ (Neumann type),
and that the precipitations on $\Gt$ is modelled as a vector field $\vp:\Gt\times(0,t_{\text{fin}})\rightarrow \mathbb{R}^2$, see Figure~\ref{fig:1}. Let $\Gf$ be the part of the boundary where the water flux is given, i.e., $\Gf=\Gl\cup \Gb$ for the considered case; 
if needed, imposing the pressure head on $\Gb$ or on another part of $\partial\Omega$ can be done using standard techniques.
As already hinted in the introduction, neither the flux nor the pressure head alone can be considered as given on $\Gt$. Indeed, the sole imposition of the flux does not exclude ponding because it does not constrain the pressure head; similarly, imposing the pressure head can lead to paradoxical results where the water influx exceeds the precipitated water.
To accommodate these considerations we introduce \[Q(\vq) = (\vp - \vq)\cdot\vn,\] that is the disparity between the precipitation rate and Darcy flux across the surface, and we impose the following conditions
\begin{equation}\label{eq:model1:3}
\psi \le 0, \quad
Q \le 0, \quad
Q\psi = 0
\end{equation}
which represent the mathematical formulation of the already-discussed seepage conditions.
Note that $Q\le 0$ allows outward flow as necessary. 
Crucially, the third condition, a complementarity requirement, ensures that at least one of the two constraints is zero. In other words, we either impose a Dirichlet boundary condition ($\psi=0$) or a Neumann boundary condition ($\vq\cdot\vn = \vp\cdot\vn$). 
These conditions bear conceptual resemblance to classical frictionless contact conditions, specifically the Karush-Kuhn-Tucker conditions encountered in contact mechanics; see, for instance,~\cite{kikuchi1988contact, wriggers2004computational, ALART1991353} for further elaboration.

\section{Discretizations}\label{sec:discretization}
In this section we describe the two already-mentioned discretizations for the seepage face condition: we will refer to the first one as ``non-hybridized'' and to the second one as ``hybridized''. Both variants are described in Section~\ref{sec:sec_discr_nitsche_part}, with details on their numerical resolution 
(specifically, non the treatment of the non-linearities that they contain) being provided in Section \ref{sec:sec_non-linearities}. 
Section~\ref{sec:sec_discr_mixed_form} contains instead preliminary material such as the discretization of the domain and the discretization of the problem with prescribed pressure head on $\Gt$, from which the two discretizations are obtained by adding the imposition of the seepage face boundary condition. 


\subsection{Mixed finite element discretization}\label{sec:sec_discr_mixed_form}

Let $\Omega_h$ denote a triangulation approximating the domain $\Omega$. The set of mesh edges in the triangulation is denoted by $\mathcal{E}$. Each element $T$ of the triangulation is assumed to conform to the subdivision of the boundary into $\Gb$, $\Gl$, $\Gt$, thereby producing the corresponding discrete boundaries $\Gbh$, $\Glh$, $\Gth$. Throughout, quantities with the subscript $h$ indicate the corresponding spatial discrete quantity.

Regarding the time discretization, we use a superscript $n$ to denote time discrete quantities, where the time corresponding to a generic $n$ is $t^n=n\Delta t$. Consequently, fully-discrete quantities are indicated with both the subscript $h$ and superscript $n$. For instance, the discrete counterpart of the pressure head $\psi$ at time $t^n$ is denoted as $\psi_h^n$.
Throughout the rest of the paper $(\cdot,\cdot)_{\Sigma}$ denotes the $L^2(\Sigma)$-scalar product.


On the triangulation $\Omega_h$ we introduce the finite element spaces $\vV_h$ and $W_h$ for approximating the Darcy flux $\vq$ and pressure head $\psi$, respectively. We choose the div-conforming inf-sup stable pair \cite{brezzifortin1991} given by the lowest order Raviart-Thomas for $\vq$ and the space of piecewise constant functions for $\psi$,
so that the water mass is conserved.
Specifically, 
\begin{equation}
    \vV_h = \{\vv\in\RT(\Omega_h): \vv\cdot\vn|_{\Gfh}=0 \}, \qquad W_h = \P(\Omega_h).
\end{equation}
Additionally, we denote by $\{\boldsymbol{\phi}_i\}_{i=1}^{\dim\vV_h}$ the set of basis functions for $\vV_h$ and by $\{\phi_i\}_{i=1}^{\dim W_h}$ the set of element-wise constant basis functions for $W_h$. 

Using an implicit Euler scheme in time, the weak form of \eqref{eq:model1:1} in the space $\vV_h$ and $W_h$ leads to the following nonlinear problem that relates $(\vq_h^{n+1},\psi_h^{n+1})$ to $(\vq_h^{n},\psi_h^{n})$: 
\begin{weakform}[Richards' equation]
Find $(\vq_h^{n+1},\psi_h^{n+1}) \in\vV_h\times W_h $ such that for all $(\vv_h, w_h)\in\vV_h\times W_h$, 
\begin{equation}\label{sec:sol:eq_1}
\begin{aligned}
a(\vq_h^{n+1},\vv_h; \psi^{n+1}_h) + b(\psi_h^{n+1},\vv_h) + e(\psi_h^{n+1},\vv_{h}) &=d(\vv_h), \\
-b(w_h,\vq_h^{n+1}) + c(w_h,\psi_h^{n+1}) &=c(w_h,\psi_h^{n}),
\end{aligned}
\end{equation}
\end{weakform}
where
\[
\begin{aligned}
a(\vq_h^{n+1},\vv_h; \psi^{n+1}_h) &= (\tK^{-1}(\psi^{n+1}_h)\vq_h^{n+1},\vv_h)_{\Omega_h}, &&& 
b(\psi_h^{n+1},\vv_h) &= -(\psi_h^{n+1},\nabla\cdot\vv_h)_{\Omega_h}, \\
c(w_h,\psi_h^{n+1}) &= \dfrac{1}{\Delta t}(\theta(\psi_h^{n+1}),w_h)_{\Omega_h},  &&& 
d(\vv_h) &= -(\nabla z,\vv_h)_{\Omega_h},
\\ e(\psi_h^{n+1},\vv_{h}) &= (\psi_h^{n+1},\vv_{h}\cdot \vn_h)_{\Gth}.
\end{aligned}
\]
If the pressure head on $\Gth$ was given, replacing $\psi_h^{n+1}$ with it in $e$ would result in the standard weak imposition of the Dirichlet data for a problem in mixed form.
Adding a similar term for $\Gb$ or $\Gl$ allows dealing with prescribed pressure head values on other parts of boundary.
The non-hybridized and the hybridized formulations differ in how the term $e(\psi_h^{n+1},\vv_{h})$ is modified to get a system that completely determines $(\vq_h^{n+1},\psi_h^{n+1})$ so that a discrete version of \eqref{eq:model1:3} is satisfied. Note that in addition to the non-linearity of the seepage boundary condition \eqref{eq:model1:3} (due to the fact that we do not know the extent of the Neumann and Dirichlet portions of $\Gth$), both methods must also cope with the nonlinear behavior of $a$ and $c$ that depend on the Van~Genuchten constitutive laws. 
 
\subsection{Numerical treatment of the seepage faces}\label{sec:sec_discr_nitsche_part}

In this section, we describe how the mixed finite element formulation introduced in Section~\ref{sec:sec_discr_mixed_form} is modified to accommodate the conditions~\eqref{eq:model1:3} at the top of the soil domain. As already mentioned, our approach is inspired by the Nitsche-based finite element scheme for contact problems developed in \cite{chouly2013nitsche, chouly2015symmetric, chouly2017overview, chouly2018unbiased}. 

We begin by noting that, as proved in \cite{chouly2013nitsche}, for any pair of real numbers $a$, $b$ the condition $a,b\le 0$ and $ab=0$ is equivalent to 
\[a=-\frac1\gamma\lbrack b - \gamma a\rbrack_+,\]
where $\gamma>0$ and $[\cdot]_{+} = \max{\left(0,\cdot\right)}$.
The seepage conditions~\eqref{eq:model1:3} on $\Gt$ are then equivalent to any of the following
\begin{align}
\label{eq:penaltyModel_1} \psi &= - \frac{1}{\gamma} \lbrack Q - \gamma \psi \rbrack_{+},
\\\label{eq:penaltyModel_2}
Q &= - \frac{1}{\gamma_{\text{hyb}}} \lbrack \psi - \gamma_{\text{hyb}} Q\rbrack_{+},
\end{align}
where $\gamma$ and $\gamma_{\text{hyb}}$ are positive functions $\Gt\rightarrow \mathbb{R}_{>0}$.
These are used in the discretized weak form setting $\gamma$ and $\gamma_{\text{hyb}}$ to edge-wise constant functions 
\[\gamma=\gamma_0 h_{e},\qquad  \gamma_{\text{hyb}} = \gamma_{0,\text{hyb}}/h_{e},\]
where $h_e$ is the edge length and $\gamma_0,\, \gamma_{0,\text{hyb}}$ are positive parameters.

The non-hybridized scheme is obtained by replacing $\psi_{h}^{n+1}$ with the expression in \eqref{eq:penaltyModel_1} for the boundary term in \eqref{sec:sol:eq_1} leading to the following
nonlinear equations that determines $\vq_h^{n+1},\psi_h^{n+1}$:  
\begin{weakform}[non-hybridized Richards' equation with seepage conditions]\label{weakform:nonhyb}
Find $\vq_h^{n+1},\psi_h^{n+1} \in\vV_h\times W_h$ such that for all $(\vv_h, w_h)\in\vV_h\times W_h$
\begin{equation}\label{discr:Nits-mixed}
\begin{aligned}
a(\vq_h^{n+1},\vv_h; \psi_h^{n+1}) + b(\psi_h^{n+1},\vv_h)  
+e( - \frac{1}{\gamma}\lbrack Q_h^{n+1}-\gamma \psi_h^{n+1}
\rbrack_{+},\vv_{h}) &= d(\vv_h), \\
-b(w_h,\vq_h^{n+1}) + c(w_h,\psi_h^{n+1}) &=c(w_h,\psi_h^{n}).
\end{aligned}
\end{equation}
\end{weakform}
\noindent The accuracy of this scheme depends on the choice of the parameter $\gamma_0$ that must be ``sufficiently small'' (we will come back to this aspect in Section \ref{sec:sec_non-linearities}).

The hybridized scheme relates instead the Dirichlet value of $\psi_h^{n+1}$ on $\Gth$
with the flux $\vq_h^{n+1}$ by introducing an auxiliary variable $\lambda_h^{n+1}\in N_h=\P(\Gth)$ (which represents the value of the pressure head on $\Gth$), and then using Equation~\eqref{eq:penaltyModel_2}. 
The basis of $N_h$ is $\{\phi_i^l\}_{i=1}^{N_b}$, with $N_b$ the number of mesh boundary edges.
The equations that determine $\vq_h^{n+1},\psi_h^{n+1}$ and $\lambda_h^{n+1}$ for the hybridized formulation are: 
\begin{weakform}[hybridized Richards' equation with seepage conditions]\label{weakform:hyb}
Find $\vq_h^{n+1},\psi_h^{n+1},\lambda_h^{n+1} \in\vV_h\times W_h\times N_h$ such that for all $(\vv_h, w_h, \eta_h)\in\vV_h\times W_h\times N_h$
\begin{equation}\label{discr:hyb}
\begin{aligned}
a(\vq_h^{n+1},\vv_h; \psi_h^{n+1}) + b(\psi_h^{n+1},\vv_h) + e(\lambda_h^{n+1}, \vv_h) &=d(\vv_h), \\
-b(w_h,\vq_h^{n+1}) + c(w_h,\psi_h^{n+1}) &= c(w_h,\psi_h^{n}), \\
e(\eta_h, \vq_h^{n+1}) +(\eta_h,- \frac{1}{\gamma_{\text{hyb}}}\lbrack \lambda_h^{n+1}  - \gamma_{\text{hyb}} Q_h^{n+1} \rbrack_{+})_{\Gth} &= e(\eta_h, \vp_h^{n+1}).
\end{aligned}
\end{equation}
\end{weakform}
The hybridization procedure adds further degrees of freedom on the discrete domain boundary
$\Gth$, but as shown later, its convergence does not depend on $\gamma_{0,\text{hyb}}$.

\subsection{Solving the non-linear equations}\label{sec:sec_non-linearities}
To determine $\psi_h^{n+1}$ and $\vq_h^{n+1}$ satisfying either
Weak~form~\ref{weakform:nonhyb} or \ref{weakform:hyb},
an iterative scheme with linear iterations is employed at each time step; 
let $k$ denote the iteration index of such iterative scheme. More specifically,
starting from the pair $(\vq_{h,0}^{n+1},\psi_{h,0}^{n+1}) = (\vq_{h}^{n},\psi_{h}^{n})$ 
a new pair $(\vq_{h,k+1}^{n+1},\psi_{h,k+1}^{n+1})$ is obtained by solving a linear system 
involving $(\vq_{h,k}^{n+1},\psi_{h,k}^{n+1})$ till a convergence criteria are achieved, i.e.
$(\psi_h^{n+1},\vq_h^{n+1})$ is the first pair in the sequence $(\vq_{h,k}^{n+1},\psi_{h,k}^{n+1})$, that satisfies the convergence criteria. 
As already hinted, there are two sources of nonlinearities: the Van~Genuchten constitutive laws and the
seepage face boundary conditions: the proposed method uses a different linearizations for each. 

\paragraph{Van~Genuchten non-linearities.}

The constitutive laws are either linearized using an $L$-scheme~\cite{slodicka2002robust, pop2004mixed, list2016study}, which is a fixed-point iteration, or the Newton method~\cite{bergamaschi1999mixed, lehmann1998comparison}.
A robust combination of both, as recently proposed in~\cite{stokke2023adaptive} is also possible.

Concerning the stability of the two nonlinear solvers, it has been shown in~\cite{list2016study} that, under the condition $L\ge\frac{1}{2}\sup_{\xi\in\mathbb{R}}{\partial_{\psi}\theta(\xi)}$, the $L$-scheme converges irrespectively of the initial guess, the time-step size and the mesh size. However,  its convergence rate is greatly influenced by the parameter $L$ that appears in Equation \eqref{eq:matrices-linear} below, see \cite{list2016study, mitra2019modified}. 
In this work, we will make use of the adaptive $L$-scheme strategy as described in Appendix A of~\cite{stokke2023adaptive} when dealing with this scheme. 

Regarding the Newton method, it usually suffers from the choice of the initial guess. For a $r$-H\"{o}lder continuous $\partial_{\psi}\theta$ function ($r\in(0,1]$), in~\cite{radu2006newton} it has been shown that Newton method exhibits $(1+r)^{\text{th}}$ order convergence if the following condition is met
\begin{equation}
    \Delta t \le C\theta_m^{\frac{2+r}{r}}h^d,
\end{equation}
where $h$ is the mesh size, $d$ is the spatial dimension, $C>0$ is a constant that depends on the problem at hand, and $\theta_m:=\inf{\partial_{\psi}\theta}\ge0$.

\paragraph{Seepage conditions non-linearities.}

The seepage face conditions are linearized by a semismooth Newton method.
Here semismooth refers to the fact that the Gateaux derivatives with respect of $\psi$ and
$\vq$ of both \eqref{eq:penaltyModel_1} and  \eqref{eq:penaltyModel_2} are discontinuous,
and equal to zero on the subset of $\Gf$ whenever both $\psi$ and $Q$ are zero.
This leads to a case by case computation, i.e. the method reduces to an active set strategy~\cite{hueber2008primal, wohlmuth2011variationally, berge2020finite}.
To give more details on this procedure, we need to write the linear
system to be solved at each iteration of the iterative algorithm. To this end,
we start by defining the terms that do not depend neither on the choice to consider
Weak form \ref{weakform:nonhyb} or \ref{weakform:hyb}, 
nor on the choice of the method employed to solve the Van~Genuchten non-linearities:
\begin{equation}\label{eq:matrices-common}
\begin{aligned}
    [\vA(\psi^{n+1}_{h,k})]_{i,j} &= a(\boldsymbol{\phi}_i,\boldsymbol{\phi}_j;\psi^{n+1}_{h,k}) = (\tK^{-1}(\psi^{n+1}_{h,k})\boldsymbol{\phi}_i,\boldsymbol{\phi}_j)_{\Omega_h},
\\  [\vB]_{i,j} &= b(\phi_i,\boldsymbol{\phi}_j)= -(\phi_i,\nabla\cdot\boldsymbol{\phi}_j)_{\Omega_h},
\\  [\vC(\psi_{h,k}^{n+1})]_{i\phantom{,j}} &= c(\phi_i;\psi_{h,k}^{n+1})-c(\phi_i;\psi_h^{n+1})=\dfrac{1}{\Delta t}(\theta(\psi_h^n) - \theta(\psi_{h,k}^{n+1}),\phi_i)_{\Omega_h},
\\  [\vD]_{i\phantom{,j}} &= -(\nabla z,\boldsymbol{\phi}_i)_{\Omega_h},
\\  [\vE]_{i,j} &= e(\phi_i^{l}, \boldsymbol{\phi}_j) =(\phi_i^{l}, \boldsymbol{\phi}_j\cdot\vn_h)_{\Gth}.
\end{aligned}
\end{equation}
The terms depending on the Van~Genuchten linearization scheme are 
\begin{equation}\label{eq:matrices-linear}
\begin{aligned}
    [\vN_{\mathcal{L}}]_{i,j} &= 
    \dfrac{1}{\Delta t}\begin{cases}
        (L \phi_i,\phi_j)_{\Omega_h}, &\text{for the $L$-scheme}, \vspace{2pt} \\
        (\partial_{\psi}\theta(\psi_{h,k}^{n+1}) \phi_i,\phi_j)_{\Omega_h}, &\text{for the Newton method},
    \end{cases}
\\      [\vB_{\mathcal{L}}]_{i,j} &= 
    \dfrac{1}{\Delta t}\begin{cases}
        [\vB]_{j,i}, &\text{for the $L$-scheme}, \vspace{2pt} \\
        [\vB]_{j,i} + (\partial_{\psi}\tK^{-1}(\psi_{h,k}^{n+1})\phi_j\vq_{h,k+1}^{n+1}\cdot\boldsymbol{\phi}_i)_{\Omega_h}, &\text{for the Newton method},
    \end{cases}
\end{aligned}
\end{equation}
where $L$ is the $L$-scheme parameter.
The blocks common to both the non-hybridized and the hybridized method are then 
the matrix $\mathbb{A}$ and the vector $\mathbb{B}$ 
\begin{equation}\label{eq:nonlin_RT0P0}
\mathbb{A} := \begin{bmatrix}
\vA(\psi^{n+1}_{h,k}) & \vB_{\mathcal{L}} \\
-\vB & \vN_{\mathcal{L}}\\
\end{bmatrix},\qquad \mathbb{B}
:=
\begin{bmatrix}
\vD \\
\vC(\psi_{h,k}^{n+1}) + \vN_{\mathcal{L}}\psi_{c,k}^{n+1}\\
\end{bmatrix}.
\end{equation}
Note that in the definition of $\mathbb{B}$ we used the symbol $\psi_{c,k}^{n+1}$ to denote
the vector of finite element degrees of freedom representing the finite element function
$\psi_{h,k}^{n+1}$; analogous notation (i.e., replacing the subscript $h$ with $c$)
will be used throughout the rest of this section with the same meaning for the other unknown fields $\vq$ and $\lambda$. 
We are now ready to discuss the treatment of the seepage boundary conditions. For convenience, 
we now split the discussion in two blocks: first the non-hybridized method (Weak form \ref{weakform:nonhyb}), 
and then the hybridized one (Weak form \ref{weakform:hyb}).

\paragraph{Non-hybridized method.}
The blocks corresponding to $e(-\frac1\gamma\lbrack Q_{h,k}^{n+1}-\gamma\psi_{h,k}^{n+1}\rbrack_+,\vv_h)$ and to its the Gateaux derivatives of  with respect to $\vq_{h,k}^{n+1}$ 
and to $\psi_{h,k}^{n+1}$ are denoted respectively by 
$\vH_{\text{no-hyb}}$, $\vF_{\text{no-hyb}}$ and $\vG_{\text{no-hyb}}$, whose expressions are:
\begin{equation}\label{eq:matrices-no-hyb}
\begin{aligned}
[\vH_{\text{no-hyb}}(\vq_{h,k}^{n+1},\psi_{h,k}^{n+1})]_{i\phantom{,j}} &= (\frac{1}{\gamma}\lbrack Q_{h,k}^{n+1}-\gamma \psi_{h,k}^{n+1}) \rbrack_{+},\boldsymbol{\phi}_i\cdot\vn_h)_{\Gth},\\
[\vF_{\text{no-hyb}}(\vq_{h,k}^{n+1},\psi_{h,k}^{n+1})]_{i,j}&=\frac{1}{\gamma}({\mathbb{1}}_{\GthN}\boldsymbol{\phi_j}\cdot \vn_h,\boldsymbol{\phi_i}\cdot\vn_h)_{\Gth},\\  [\vG_{\text{no-hyb}}(\vq_{h,k}^{n+1},\psi_{h,k}^{n+1})]_{i,j}&=
({\mathbb{1}}_{\GthN} \phi_j,\boldsymbol{\phi_i}\cdot\vn_h)_{\Gth},
\end{aligned}
\end{equation}
where $\mathbb{1}_{\GthN}$ is the characteristic function of 
$\GthN =\{Q_{h,k}^{n+1}-\gamma \psi_{h,k}^{n+1} \ge 0\}\subseteq \Gth$.
Notice that 
$\psi_{h,k}^{n+1}$ is piecewise constant so its value on $\Gth$ (its trace) can be chosen as its value in the middle of the mesh elements touching $\Gth$,
for any $\vq_h^n \in V_h$ the corresponding quantity $Q_h^n=(\vp_h^n-\vq_h^n)\cdot \vn_h$ is edge-wise constant on $\Gth$ and consequently it is in $\P(\Gfh)$, so that
$\vH_{\text{no-hyb}}$ does not have projection errors.
We have now all terms that define the non-hybridized iteration:
\begin{equation}\label{discr_Nitsche_mixed}
\begin{aligned}
&\mathbb{A}_{\text{no-hyb}}\begin{bmatrix}
\vq^{n+1}_{c,k+1}\vspace{3pt} \\
\psi_{c,k+1}^{n+1}\\
\end{bmatrix} = \mathbb{B}_{\text{no-hyb}}, \\
&\mathbb{A}_{\text{no-hyb}} := \mathbb{A} + \begin{bmatrix}
\vF_{\text{no-hyb}}(\vq_{h,k}^{n+1},\psi_{h,k}^{n+1}) & \vG_{\text{no-hyb}}(\vq_{h,k}^{n+1},\psi_{h,k}^{n+1}) \\
\mathbf{0} & \mathbf{0} \\
\end{bmatrix}, \\
&\mathbb{B}_{\text{no-hyb}} := 
\mathbb{B} + 
\begin{bmatrix}
\vF_{\text{no-hyb}}(\vq_{h,k}^{n+1},\psi_{h,k}^{n+1})\vq_{c,k}^{n+1} + \vG_{\text{no-hyb}}(\vq_{h,k}^{n+1},\psi_{h,k}^{n+1})\psi_{c,k}^{n+1} + \vH_{\text{no-hyb}}(\vq_{h,k}^{n+1},\psi_{h,k}^{n+1})\\
\mathbf{0}
\end{bmatrix}.
\end{aligned}
\end{equation}
For edges contained in $\GthN$ we obtain the following equation, where $\gamma$ acts as a ``penality'' parameter for the weak imposition of the homogeneous Neumann boundary condition and thus influences the accuracy of the solution:
\begin{equation}\label{eq:discr_weak_constr}
\begin{aligned}
    \vA(\psi_{h,k}^{n+1})\vq_{c,k+1}^{n+1} + \dfrac{1}{\gamma}\vq_{h,k+1}^{n+1}\cdot\vn_h &=\dfrac{1}{\gamma}\vp^{n+1}_h\cdot\vn_h.
\end{aligned}
\end{equation}
If an edge is instead not contained in $\GthN$, the corresponding entries in $\vF_{\text{no-hyb}}(\vq_{h,k}^{n+1},\psi_{h,k}^{n+1}),\vG_{\text{no-hyb}}(\vq_{h,k}^{n+1},\psi_{h,k}^{n+1})$, $\vH_{\text{no-hyb}}(\vq_{h,k}^{n+1},\psi_{h,k}^{n+1})$ are null and the system above for the $(k+1)$-th iterate reduces to the weak imposition of the homogeneous Dirichlet boundary condition on $\Gth$.

\paragraph{Hybridized method.}
The method with the hybridization procedure on the boundary $\Gt$, i.e., the nonlinear set of equation~\eqref{discr:hyb}
contains the blocks corresponding to $(\eta_h,- \frac{1}{\gamma_{\text{hyb}}}\lbrack \lambda_h^{n+1}  - \gamma_{\text{hyb}} Q_h^{n+1} \rbrack_{+})_{\Gth}$ and to its derivatives with respect to $\vq_{h,k}^{n+1}$ 
 and to $\psi_{h,k}^{n+1}$, which are denoted respectively by 
$\vH_{\text{hyb}}$, $\vF_{\text{hyb}}$ and $\vG_{\text{hyb}}$, and whose expressions are: 
\begin{equation}\label{eq:matrices-hyb}
\begin{aligned}
[\vH_{\text{hyb}}(\vq_{h,k}^{n+1},\lambda_{h,k}^{n+1})]_{i\phantom{,j}} &= (\frac{1}{\gamma_{\text{hyb}}}\lbrack \lambda_{h,k}^{n+1}-\gamma_{\text{hyb}} Q_{h,k}^{n+1}\rbrack_{+},\phi_i^{l})_{\Gth},\\
[\vF_{\text{hyb}}(\vq_{h,k}^{n+1},\lambda_{h,k}^{n+1})]_{i,j}&=-({\mathbb{1}}_{\GthH}\boldsymbol{\phi_j}\cdot \vn_h,\boldsymbol{\phi_i}\cdot\vn_h)_{\Gth},\\
[\vG_{\text{hyb}}(\vq_{h,k}^{n+1},\lambda_{h,k}^{n+1})]_{i,j}&= -\frac{1}{\gamma_{\text{hyb}}}({\mathbb{1}}_{\GthH} \phi_j,\boldsymbol{\phi_i}\cdot\vn_h)_{\Gth},
\end{aligned}
\end{equation}
where $\mathbb{1}_{\GthH}$ is the characteristic function of $\GthH =\{\lambda_{h,k}^{n+1}-\gamma_{\text{hyb}}Q_{h,k}^{n+1}\ge 0 \}\subseteq \Gth$.
Its linear iteration is defined by
\begin{equation}\label{discr:lin_sys_new_hyb}
\begin{aligned}
& \mathbb{A}_{\text{hyb}}\begin{bmatrix}
\vq^{n+1}_{c,k+1}\vspace{3pt} \\
\psi_{c,k+1}^{n+1}\vspace{3pt}\\
\lambda_{c,k+1}^{n+1}\\
\end{bmatrix} = \mathbb{B}_{\text{hyb}}, \\
&\mathbb{A}_{\text{hyb}} := 
\begin{bNiceArray}{cw{c}{.1cm}c|c}[margin]
\Block{2-2}<\Large>{\mathbb{A}} & & &\mathbf{0} \\
& &  &\mathbf{0} \\
\hline
\mathbf{0} & & \mathbf{0} & \mathbf{0}
\end{bNiceArray}
+
\begin{bNiceArray}{cw{c}{.1cm}c|c}[margin]
\Block{2-2}<\Large>{\mathbf{0}} & & &\vE^T \\
& &  &\mathbf{0} \\
\hline
\vE+\vF_{\text{hyb}}(\vq_{h,k}^{n+1},\lambda_{h,k}^{n+1}) & & \mathbf{0} & \vG_{\text{hyb}}(\vq_{h,k}^{n+1},\lambda_{h,k}^{n+1})
\end{bNiceArray}, \\
&\mathbb{B}_{\text{hyb}} := 
\begin{bmatrix}
\mathbb{B} \\
\mathbf{0}
\end{bmatrix} + 
\begin{bmatrix}
\mathbf{0} \\
\vE\vp_{c}^{n+1} + \vF_{\text{hyb}}(\vq_{h,k}^{n+1},\lambda_{h,k}^{n+1})\vq_{c,k}^{n+1} + \vG_{\text{hyb}}(\vq_{h,k}^{n+1},\lambda_{h,k}^{n+1})\lambda_{c,k}^{n+1} + \vH_{\text{hyb}}(\vq_{h,k}^{n+1},\lambda_{h,k}^{n+1})
\end{bmatrix}
.
\end{aligned}
\end{equation}

In this case, on the faces belonging $\GthH$ the iteration scheme imposes
the homogeneous Neumann boundary conditions strongly and the homogeneous Dirichlet weakly on the proper subset of edges.
Indeed, for an edge not contained in $\GthH$ the quantities $\vF_{\text{hyb}}(\vq_{h,k}^{n+1},\lambda_{h,k}^{n+1}),\vG_{\text{hyb}}(\vq_{h,k}^{n+1},\lambda_{h,k}^{n+1})$, $\vH_{\text{hyb}}(\vq_{h,k}^{n+1},\lambda_{h,k}^{n+1})$ nullify, and the third row reads
\begin{equation}\label{discr_cond_1}
\begin{aligned}
    \vq_{h,k+1}^{n+1}\cdot\vn_h = \vp_{h}^{n+1}\cdot\vn_h, \quad \text{on } \Gth\setminus\GthH,
\end{aligned}
\end{equation}
whereas if a mesh edge is contained in $\GthH$ the system imposes the Dirichlet condition
\begin{equation}\label{discr_cond_2}
    \lambda_{h,k+1}^{n+1} = 0, \quad \text{on } \GthH.
\end{equation}
Differently from what happens for the discretization~\eqref{discr_Nitsche_mixed}, both~\eqref{discr_cond_1} and~\eqref{discr_cond_2} do not depend on the penalty parameter $\gamma_{0,\text{hyb}}$. This is a great advantage of the hybrid method because no penalty parameter have to be calibrated and any $\gamma_{0,\text{hyb}}>0$
guarantees the proper imposition of the seepage face conditions.
The stability of the method with respect of $\gamma_{0,\text{hyb}}$ is verified by
the numerical examples in the next section.

\section{Numerical examples}\label{sec:examples}


The two methods described above have been implemented in the open-source Python toolbox PyGeoN~\cite{PyGeoN}, which is built on top of the Python library PorePy~\cite{keilegavlen2021porepy}. PorePy itself utilizes Gmsh~\cite{geuzaine2009gmsh} for meshing procedures. 
The sparse linear systems~\eqref{discr_Nitsche_mixed} and~\eqref{discr:lin_sys_new_hyb} that arise at each nonlinear iteration are solved using the direct solver implemented in SciPy's \textit{spsolve} routine~\cite{2020SciPy-NMeth}. \\
The stopping criteria for the nonlinear solver is 
\begin{equation}\label{eq:stopping_criterion}
    \eta_{\text{lin}}^{k+1} := |||{\psi}_{h,k+1}^{n+1} - {\psi}_{h,k}^{n+1}|||_{\mathcal{L}} \leq \epsilon_{A}
\end{equation}
where $\epsilon_{A}$ is the absolute error tolerance (if not otherwise stated we will take $\epsilon_A=10^{-7}$ in the simulations we present below), and $|||\cdot|||_{\mathcal{L}}$ is the particular $H^1(\Omega_h)$ equivalent-norm associated with the scheme chosen to solve the-linearities, i.e., either the $L$-scheme or the Newton method as discussed in Section \ref{sec:sec_non-linearities}.
Such norm is defined as follows:
\begin{equation}\label{eq:eta_error}
\begin{aligned}
    \eta_{\text{lin}}^{k+1} &= \left(\int_{\Omega_h} \mathbb{L}({\psi}_{h,k+1}^{n+1} - {\psi}_{h,k}^{n+1})^2 + \Delta t \left|\sqrt{K({\psi}_{h,k}^{n+1})} (\tK^{-1}({\psi}_{h,k+1}^{n+1})\vq_{h,k+1}^{n+1} - \tK^{-1}({\psi}_{h,k}^{n+1})\vq_{h,k}^{n+1}) \right|^2 \right)^{\frac{1}{2}}, \\
     \mathbb{L} &= \begin{cases}
       L, & \text{for the $L$-scheme},\\
       \partial_{\psi}\theta({\psi}_{h,k}^{n+1}), & \text{for the Newton method}.
    \end{cases}
\end{aligned}
\end{equation}
Note that the norm above encapsulates the entirety of the linearization error, as shown in Section $5$ of~\cite{mitra2023guaranteed}.

\begin{figure}[h!]
\centering
\includegraphics[width=0.5\textwidth]{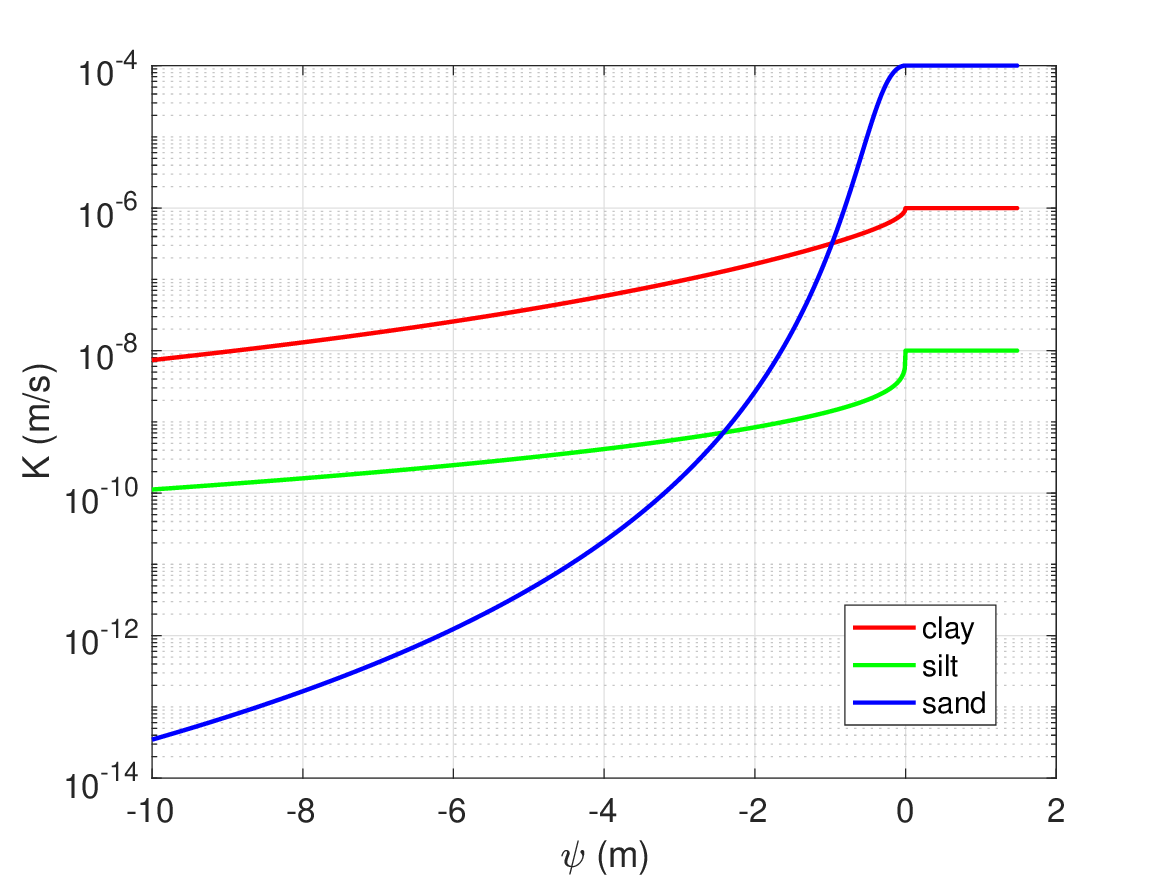}\hspace*{-.4cm}
\includegraphics[width=0.5\textwidth]{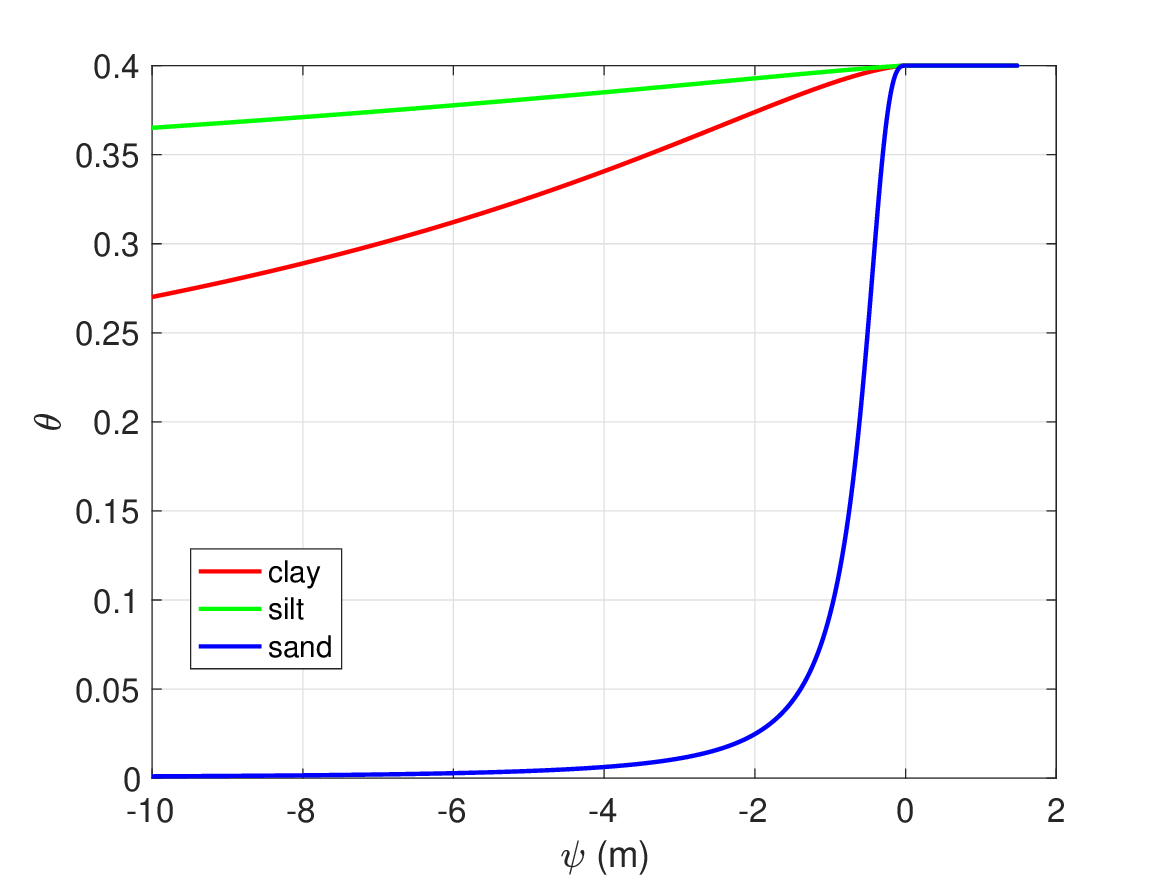}
\caption{Left panel shows of the permeability curve $K(\psi)$ for the various soil types. Right panel reports of the water content curve $\theta(\psi)$ for the various soil types. }
\label{fig:14}
\end{figure}

Various domain complexities will be considered to test the performance of the proposed discretization schemes to deal with the seepage conditions on $\Gt$. 
In Section~\ref{sec:rec_domain_1} we first consider the simplest case, which corresponds to the domain in Figure \ref{fig:1} with horizontal top surface (i.e., a rectangular domain). 
Then, in Section~\ref{sec:ideal_slope}, we consider a more realistic scenario
consisting of two plateaus with different height, connected by an inclined surface (\emph{artificial slope}).
Finally, in Section~\ref{sec:real_slope_DTM}, we show the ability of the proposed schemes to deal with natural slopes taken from DTM data. 
Regardless of the geometry of the domain, the precipitation will be assumed to be aligned with the vertical direction, i.e., $\vp=-p\mathbf{e}_z$, where $p>0$ is the actual rainfall that we assume constant and uniform.

In all cases, the domain $\Omega$ is considered filled just by one soil type, that in the proposed tests will be one of clay, slit, or sand, that are characterized by the following Van~Genuchten parameters:
\begin{center}
\begin{tabular}{cccccc}
    \textbf{Material}   & $\theta_S$ & $\theta_R$   & $\alpha$      & $m$   & $K_S$ \\
    clay                & $0.4$      &  $0.04$      & $0.2$m$^{-1}$ & $1.5$ & $10^{-6}$m/s \\  
    silt                & $0.4$      &  $0.08$      & $0.1$m$^{-1}$ & $1.2$ & $10^{-8}$m/s \\  
    sand                & $0.4$      &  $0\phantom{.00}$         & $2\phantom{.0}$m$^{-1}$   & $3\phantom{.0}$   & $10^{-4}$m/s \\      
\end{tabular}
\end{center}

In Figure~\ref{fig:14} we plot the permeability and water content curve for the various soil types considered.
The simulations are performed with a MacBook Air laptop (M1, 2020) with $8$GB of RAM.


\subsection{Rectangular domain}\label{sec:rec_domain_1}

Let us consider a simple rectangular domain $\Omega=(0,0.1)\times(0,5)$m$^2$, and set uniform and constant rainfall conditions on $\Gt$ such that $p$ does not depend on space and time. 
On the bottom $\Gb$ we apply homogeneous Dirichlet boundary conditions, and on the lateral boundaries $\Gl$ we apply no-flux boundary conditions. 
The initial condition is set equal to hydrostatic ranging from $0$m to the bottom to the value of $-5$m at the top. 
The material filling $\Omega$ is assumed to be clayey soil. 
Concerning the numerical discretization, we choose a uniform mesh size equal to $0.05$m. We set a maximum number of iterations of the nonlinear solver equal to $200$, that is never reached for the simulations in this subsection. To solve the Van Genuchten nonlinearities we use the $L$-scheme with a tolerance on the stopping criteria $\epsilon_A=10^{-5}$, 
cf. Equation \eqref{eq:stopping_criterion}. 

In the following we consider various rainfall scenarios by changing the value of $p/K_S$, i.e., of the ratio between rainfall and saturated permeability. 
In particular, if $p/K_S<1$ the rainfall is not ``strong enough'' to saturate the soil at the surface, which means that we expect the pressure head to remain negative on $\Gt$: therefore a simple Neumann condition on $\Gt$ would
be able to simulate correctly this scenario and can be used as a reference solution
to which our approaches should be compared (i.e., a consistency test). 
Conversely, if $p/K_S > 1$ the rainfall exceeds the absorption capacity of the soil, 
such that runoff will happen and the soil will saturate starting from $\Gt$. 
Our approaches are expected to correctly capture this physics, whereas 
a Neumann boundary condition would incorrectly predict a positive pressure head on $\Gt$,
indicating ponding. 

\begin{figure}[h!]
\centering
\includegraphics[width=0.5\textwidth]{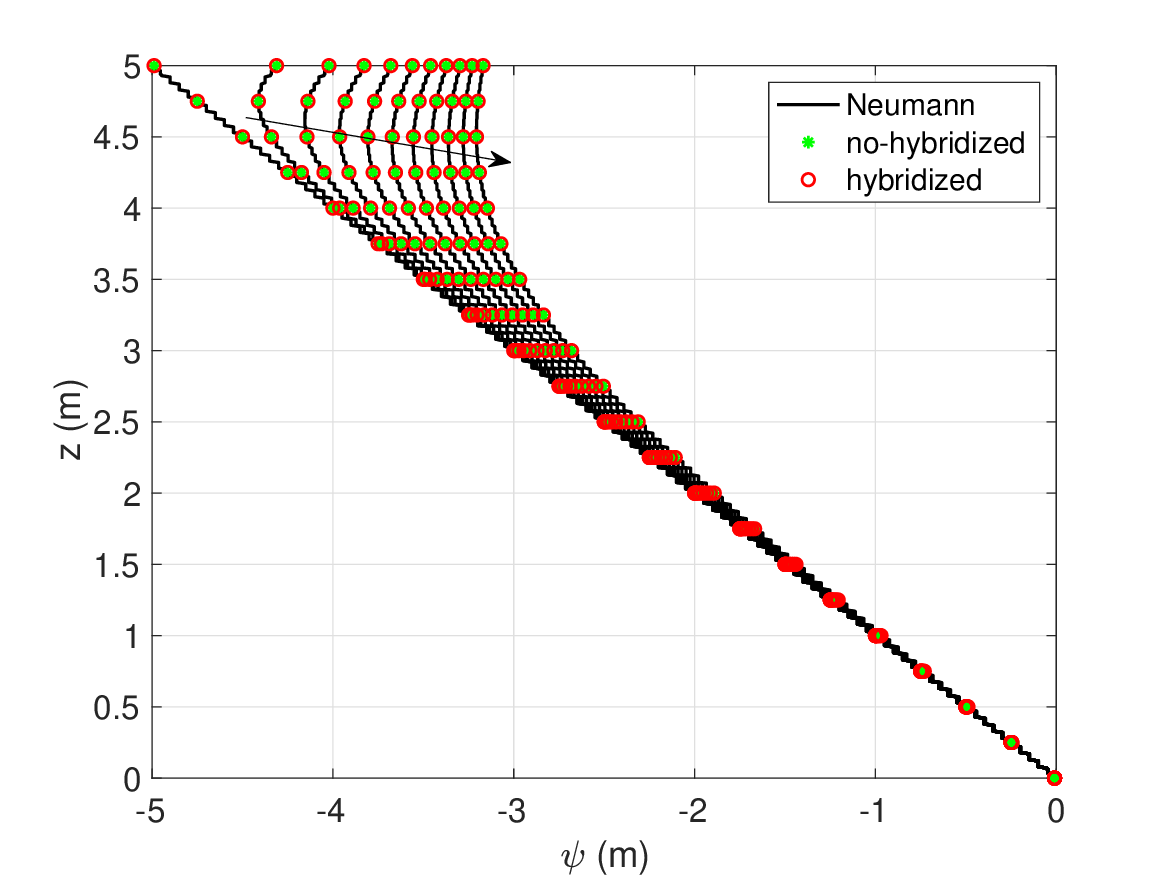}\hspace*{-.4cm}
\includegraphics[width=0.5\textwidth]{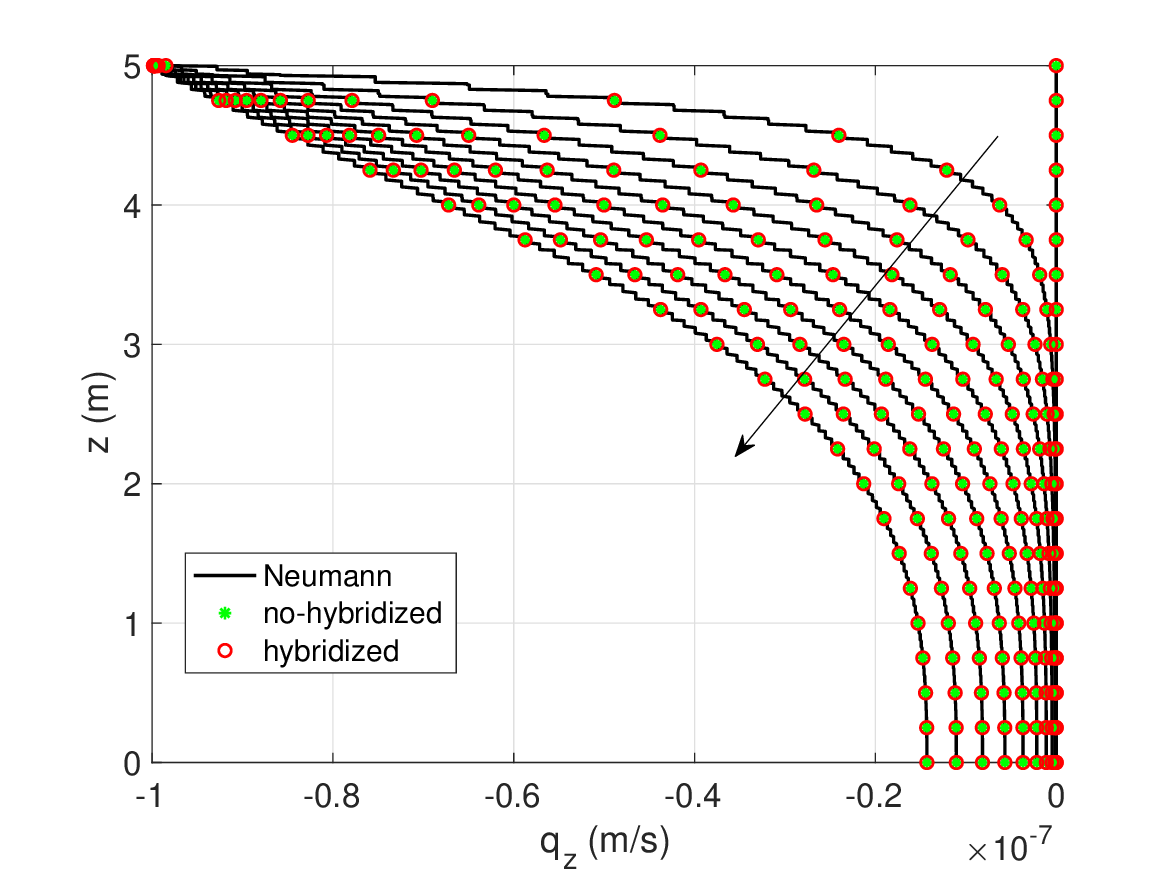}\vfill
\includegraphics[width=0.5\textwidth]{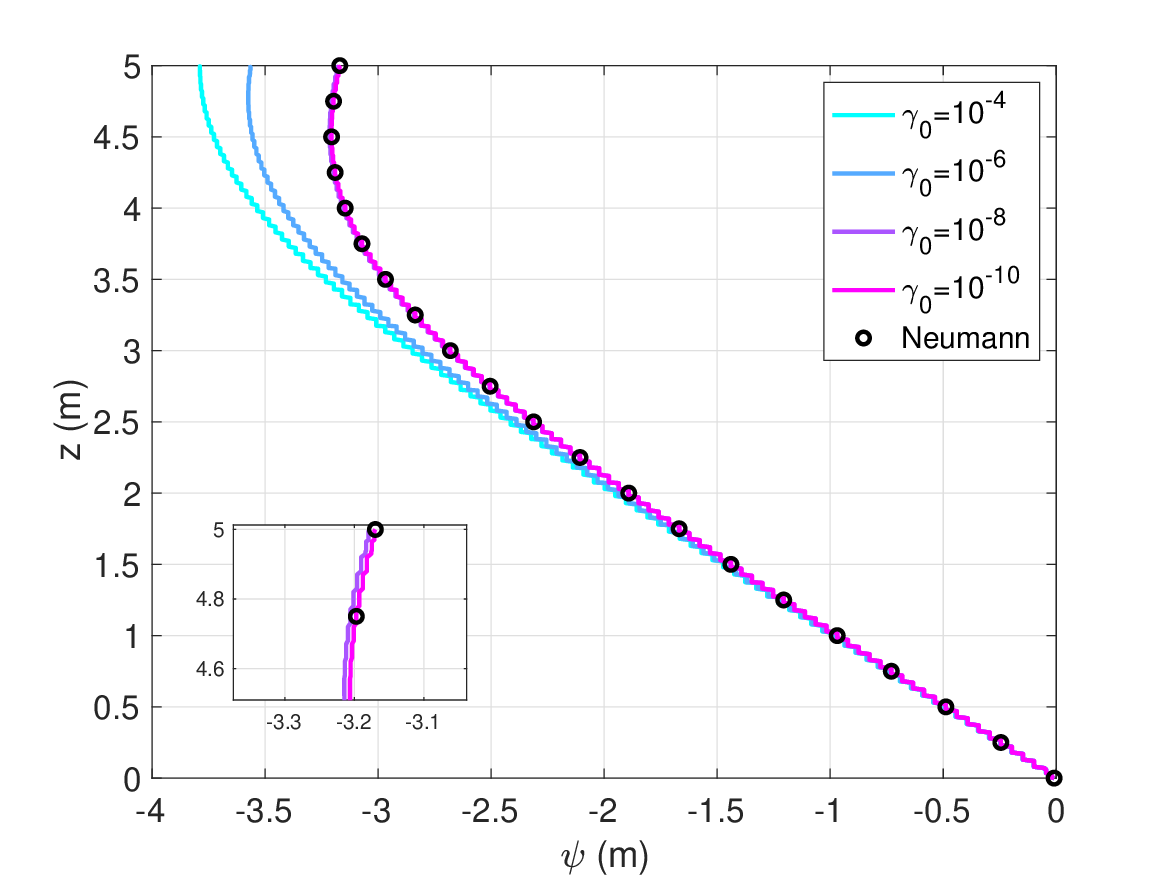}
\caption{Rectangular domain. Test with $p/K_S=0.1$. Top row: comparison between the Neumann, no-hybridized and hybridized scheme's solutions at time instants $0$h, $10$h, $20$h, \dots, $90$h, $100$h. Bottom row: sensitivity analysis on the numerical parameter $\gamma_0$ for the no-hybridized scheme on the simulation's result at final time.}
\label{fig:2}
\end{figure}

\paragraph{Case $p/K_S<1$.} The first test we carry out is characterized by a rainfall-to-saturated permeability ratio smaller than one, i.e., $p/K_S=0.1$. 
The final simulation time is set equal to $t_{\text{fin}}=100$h, and we divide the simulation in $10$ time steps, i.e. $\Delta t = 10$h. 
In Figure~\ref{fig:2} we report the results of the numerical solution extracted along the vertical line at $x=0.05$m; 
here and in the following analogous plots, the arrows indicate the trend of the pressure/flux profiles  as time increases.
In the top row we compare the results of the simulation with the non-hybrized and hybridized schemes' results at various time instants with the solution obtained imposing Neumann 
boundary conditions.
Here and in the next sets of plots, the quantity $q_z$ is the vertical component of the water flux. 
\noindent 
The results shown are obtained by setting the parameter $\gamma_0=10^{-10}$ for the non-hybridized scheme: the importance of choosing a very small value of $\gamma_0$ is evident in the panel in the bottom row, where we present a sensitivity analysis to 
$\gamma_0$ of the numerical solution at $t=t_{\text{fin}}$.

Conversely, 
we do not report any sensitivity analysis on $\gamma_{0,{\text{hyb}}}$ appearing in the hybridized scheme, that here is set to $\gamma_{0,{\text{hyb}}}=1$.
Indeed, as shown in Section~\ref{sec:sec_discr_nitsche_part},
the solution of the discrete
system of the hybridized method does not depend on it. 
We will nonetheless give numerical evidence supporting this statement later on, in
the more relevant case $p/K_S>1$.
The absolute error computed in $L^{\infty}(\Omega)$-norm with respect to the Neumann solution is approximately $10^{-3}$ for the no-hybridized scheme with $\gamma_0=10^{-10}$, while $10^{-4}$ for the hybridized one.

\paragraph{Case $p/K_S=1$.}  
This case is a limit case, in which the rainfall saturates the soil without
runoff, and the saturation front proceeds downward from $\Gt$ with the rate
given by the incoming rainfall flux.
In this case we can still compare the results of our approaches to the solution obtained
with a Neumann boundary condition; note however that in this case, contrary to the previous one,
the switch between Neumann and Dirichlet is expected to happen in our approach.
\begin{figure}[h!]
\centering
\includegraphics[width=0.5\textwidth]{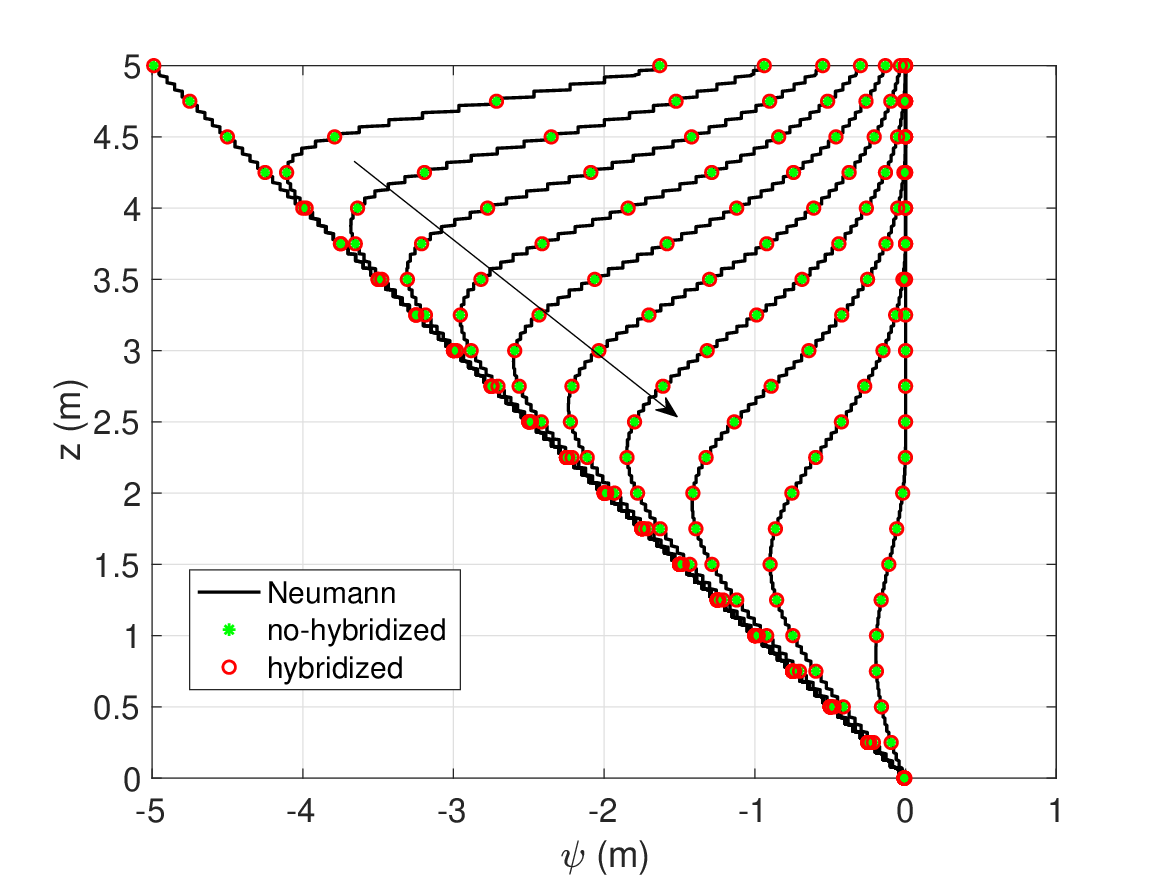}\hspace*{-.4cm}
\includegraphics[width=0.5\textwidth]{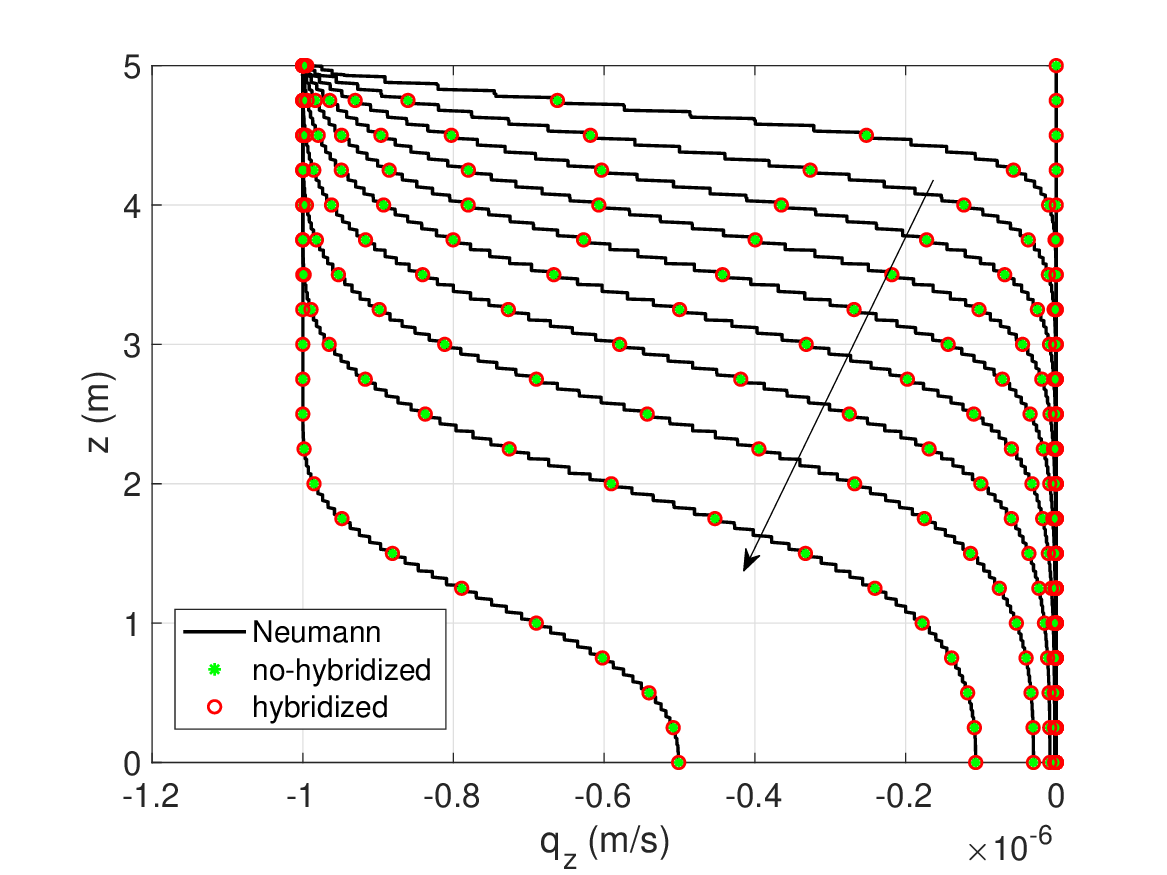}\vfill
\includegraphics[width=0.5\textwidth]{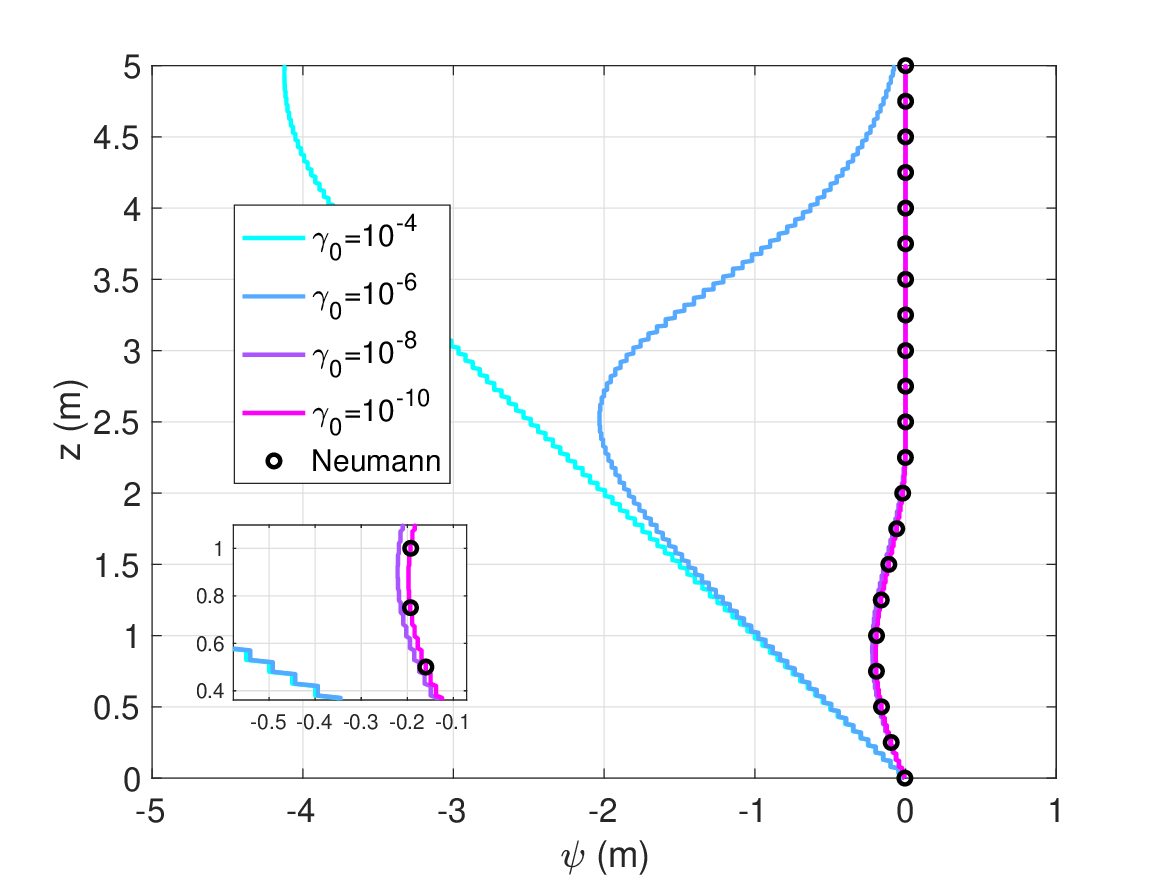}
\caption{Rectangular domain. Test with $p/K_S=1$. Top row: comparison between the Neumann, no-hybridized and hybridized scheme's solutions at time instants $0$h, $5$h, $10$h, \dots, $45$h, $50$h. Bottom row: sensitivity analysis on the numerical parameter $\gamma_0$ for the no-hybridized scheme on the simulation's result at final time.}
\label{fig:6}
\end{figure}
The final simulation time is set equal to $t_{\text{fin}}=50$h, and we take a time step $\Delta t=5$h, and the corresponding results are reported in Figure~\ref{fig:6}. 
Again, the panel in the bottom row reports a sensitivity analysis that shows the 
importance of $\gamma_0$; the results in the top row are obtained with $\gamma_0=10^{-10}$. 
\begin{figure}[h!]
\centering
\includegraphics[width=0.5\textwidth]{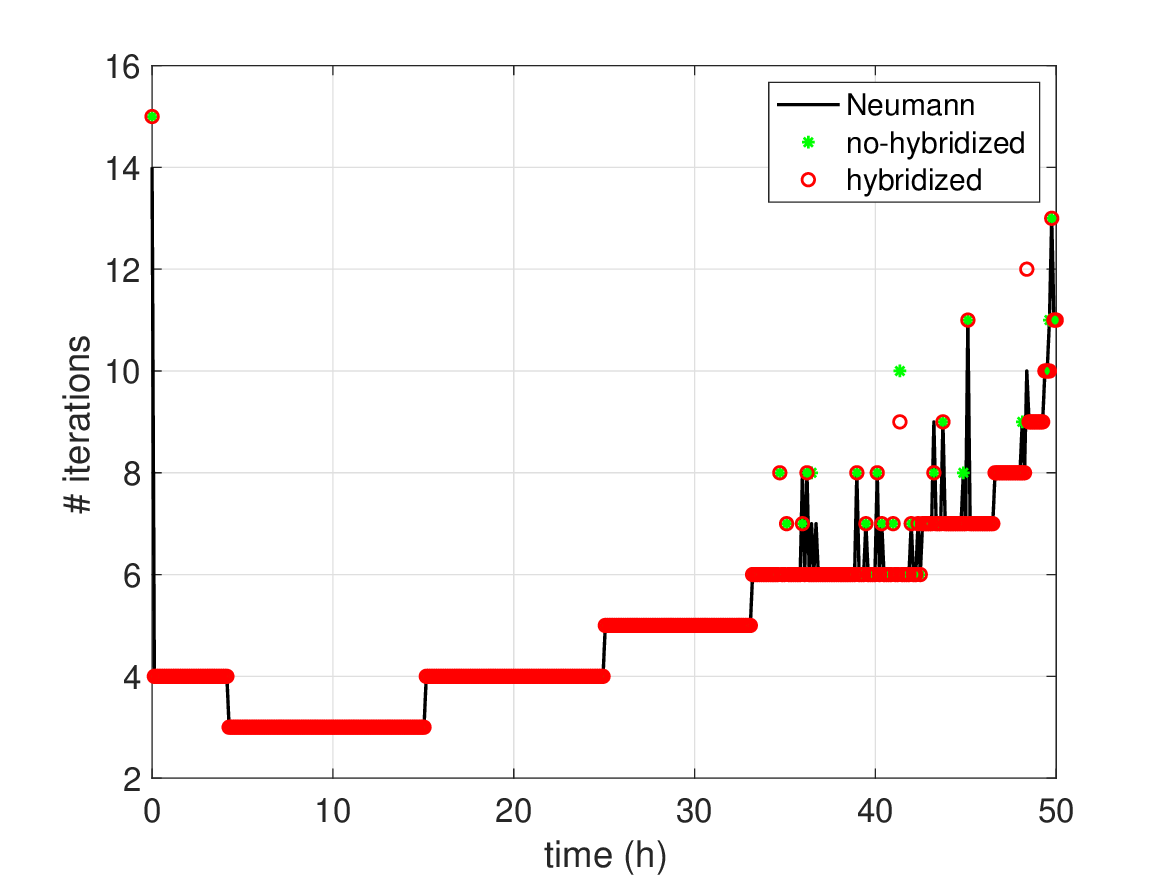}\hspace*{-.4cm}
\includegraphics[width=0.5\textwidth]{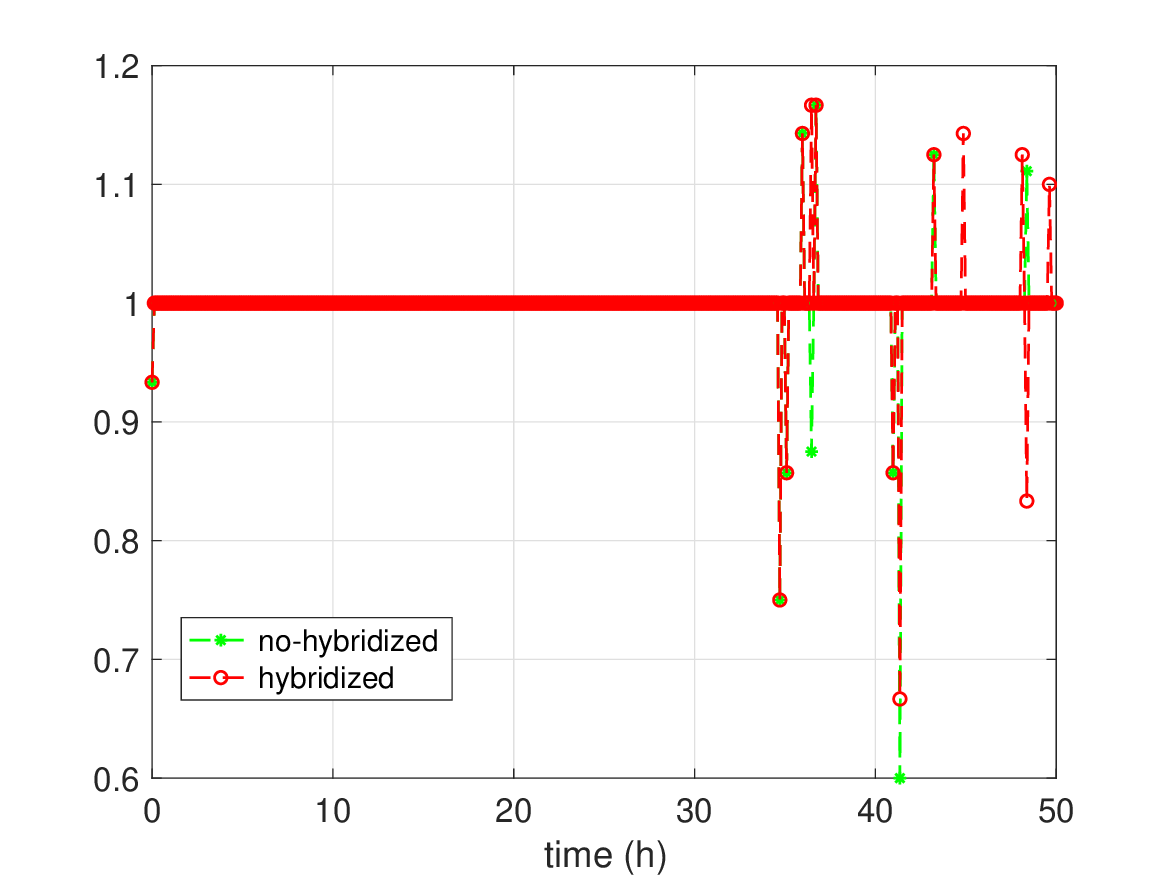}
\caption{Rectangular domain. Test with $p/K_S=1$. Left panel: number of iterations. Right panel: ratio between the number of nonlinear iterations performed by the Neumann approximation over the number required by the no-hybridized and hybridized scheme. }
\label{fig:66}
\end{figure}

The two proposed methods, the non-hybridized and hybridized one, can accurately reproduce the expected solution given by the imposition of the Neumann boundary condition on $\Gt$. The absolute error computed in $L^{\infty}(\Omega)$-norm with respect to the Neumann solution is approximately $2\cdot10^{-3}$ for the no-hybridized scheme with $\gamma_0=10^{-10}$, while $6\cdot10^{-4}$ for the hybridized one. Again, for the hybridized scheme the convergence to the Neumann solution does not depend on the choice of $\gamma_{0,{\text{hyb}}}$ so that we take it equal to $1$.

Finally, in Figure~\ref{fig:66} we show the number of nonlinear iterations during the simulation time for the three simulations performed on the left panel, and the ratio between the number of nonlinear iterations performed when considering Neumann boundary conditions over the number required by the no-hybridized and hybridized scheme on the right panel. This shows that the number of iterations required by the three methods are comparable.

\begin{figure}[h!]
\centering
\includegraphics[width=0.5\textwidth]{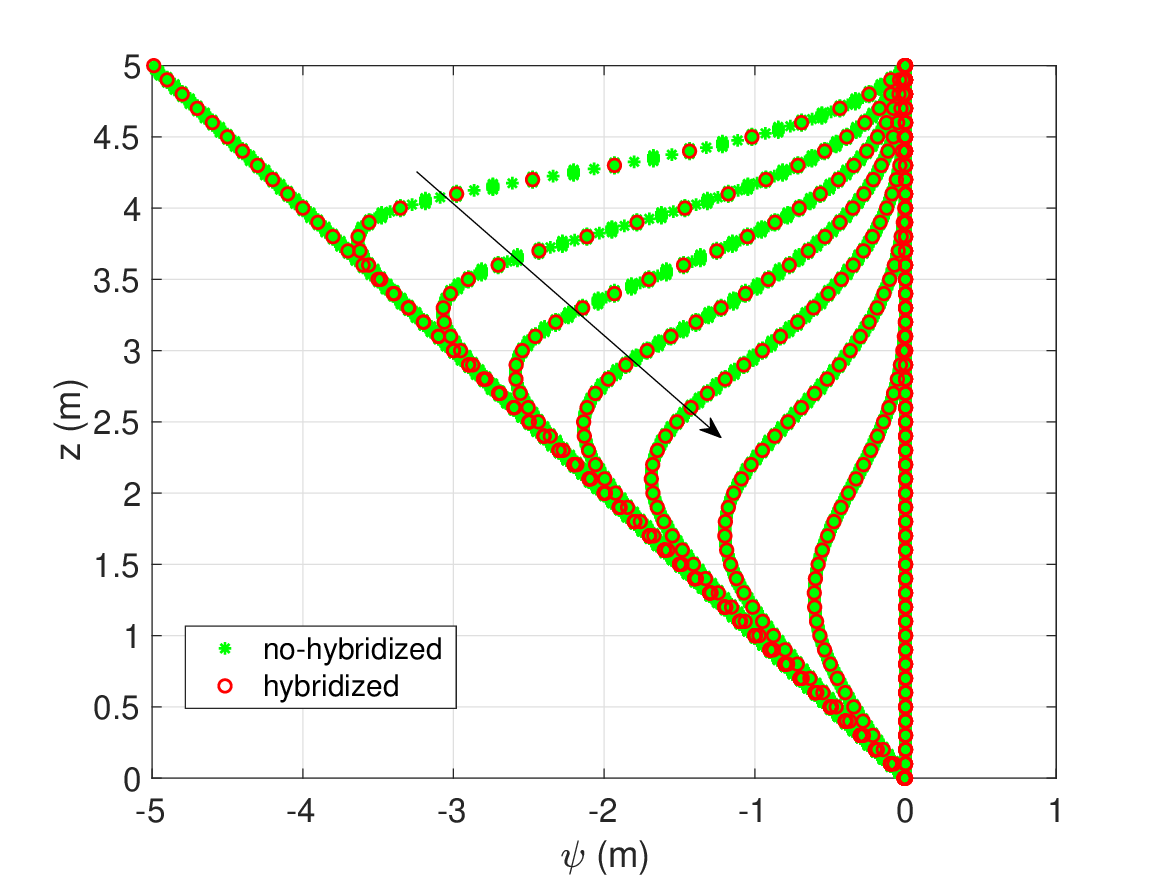}\hspace*{-.4cm}
\includegraphics[width=0.5\textwidth]{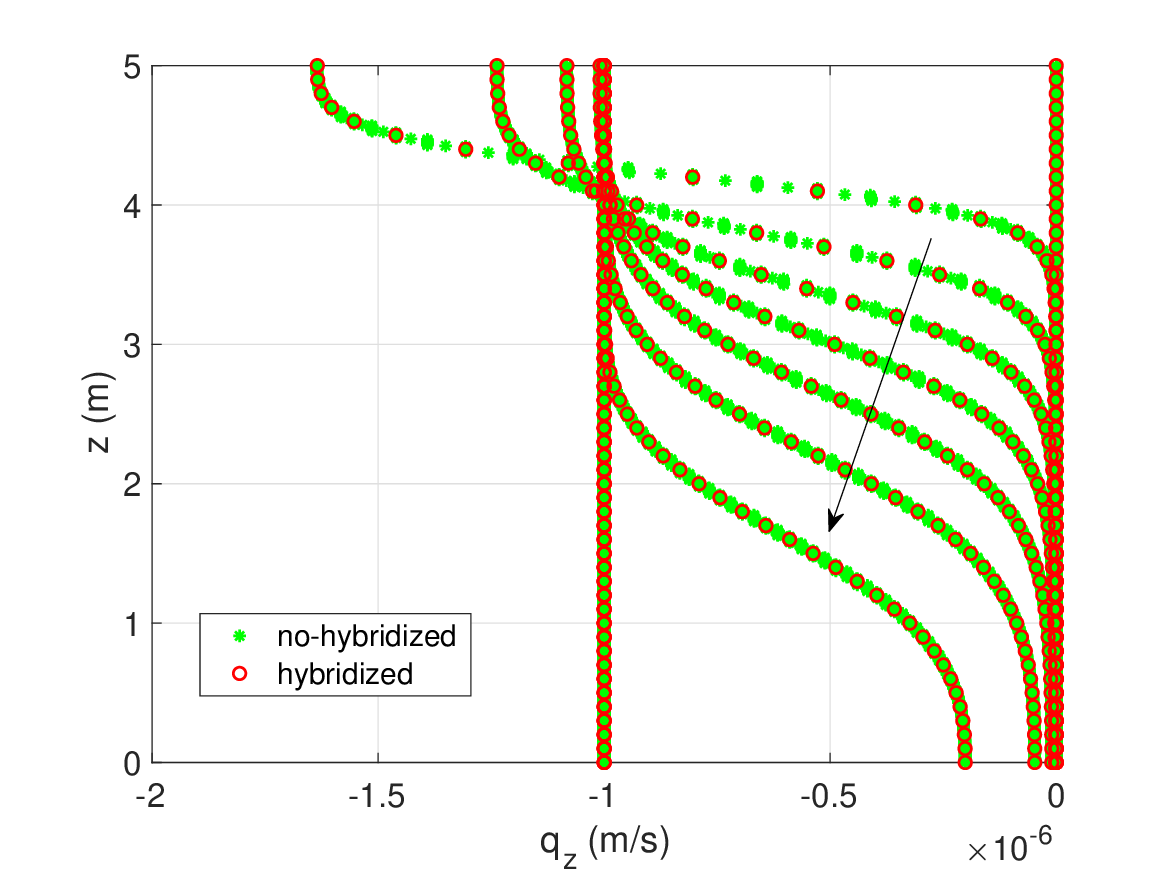}
\caption{Rectangular domain. Test with $p/K_S=10$. Comparison between the Neumann, no-hybridized and hybridized scheme's solutions at time instants $0$h, $5$h, $10$h, \dots, $45$h, $50$h. }
\label{fig:5}
\end{figure}

\paragraph{Case $p/K_S>1$.}  In the last test of this sequence, we consider a rainfall-to-saturated-permeability ratio equal to $p/K_S=10$. 
Again, the final simulation time is set equal to $t_{\text{fin}}=50$h, and we take a time step $\Delta t=5$h. In Figure~\ref{fig:5} we report the results of our analysis. We have considered $\gamma_0=10^{-10}$ and $\gamma_{0,{\text{hyb}}}=1$.
The imposition of the non-positivity of the pressure head on $\Gt$ is clearly effective. While the simulation evolves the pressure head profile approaches zero, and furthermore starting from time $t=40$h the domain $\Omega$ is completely filled by water as the pressure head becomes numerically zero everywhere.

In this third case, which is remarkably the most important from an application viewpoint, we show a table concerning the sensitivity analysis varying the parameter $\gamma_{0,{\text{hyb}}}$. This is to corroborate the theoretical result found at the end of Section~\ref{sec:sec_discr_nitsche_part} concerning the independence of the approximate solution with respect to the choice of $\gamma_{0,{\text{hyb}}}$. In Table~\ref{table}, we report the difference in $L^{\infty}(\Omega)$-norm between the solution computed with the no-hybridized scheme with $\gamma_0=10^{-10}$ and three approximations given by the hybridized method for three different parameters $\gamma_{0,{\text{hyb}}}$. As the reader can notice, the three errors are of the same order of magnitude thus not varying significantly among each other.

\begin{table}[h!]
\centering
\begin{tabular}{c | c} 
 $\gamma_{0,\text{hyb}}$ & error \\ 
 \hline
 1 & $1\cdot10^{-4}$ \\
 10$^{10}$ & $1\cdot10^{-4}$ \\
 10$^{-10}$ & $2\cdot10^{-4}$
\end{tabular}
\caption{Rectangular domain. Test with $p/K_S=10$. The error indicates the difference in $L^{\infty}(\Omega)$-norm between the no-hybridized (with $\gamma_0=10^{-10}$) and hybridized approximation for various values of the Nitsche parameter $\gamma_{0,{\text{hyb}}}$.}
\label{table}
\end{table}

\subsubsection{Relaxation of the seepage conditions for certain soil types}\label{sec:rel:soil_type}

When moving from clayey to silty soil, the temporal evolution of the pressure and flux profiles
shows some temporal discrepancies among the different discretization schemes, that can be fixed by relaxing the seepage conditions. 
\begin{figure}[h!]
\centering
\includegraphics[width=0.5\textwidth]{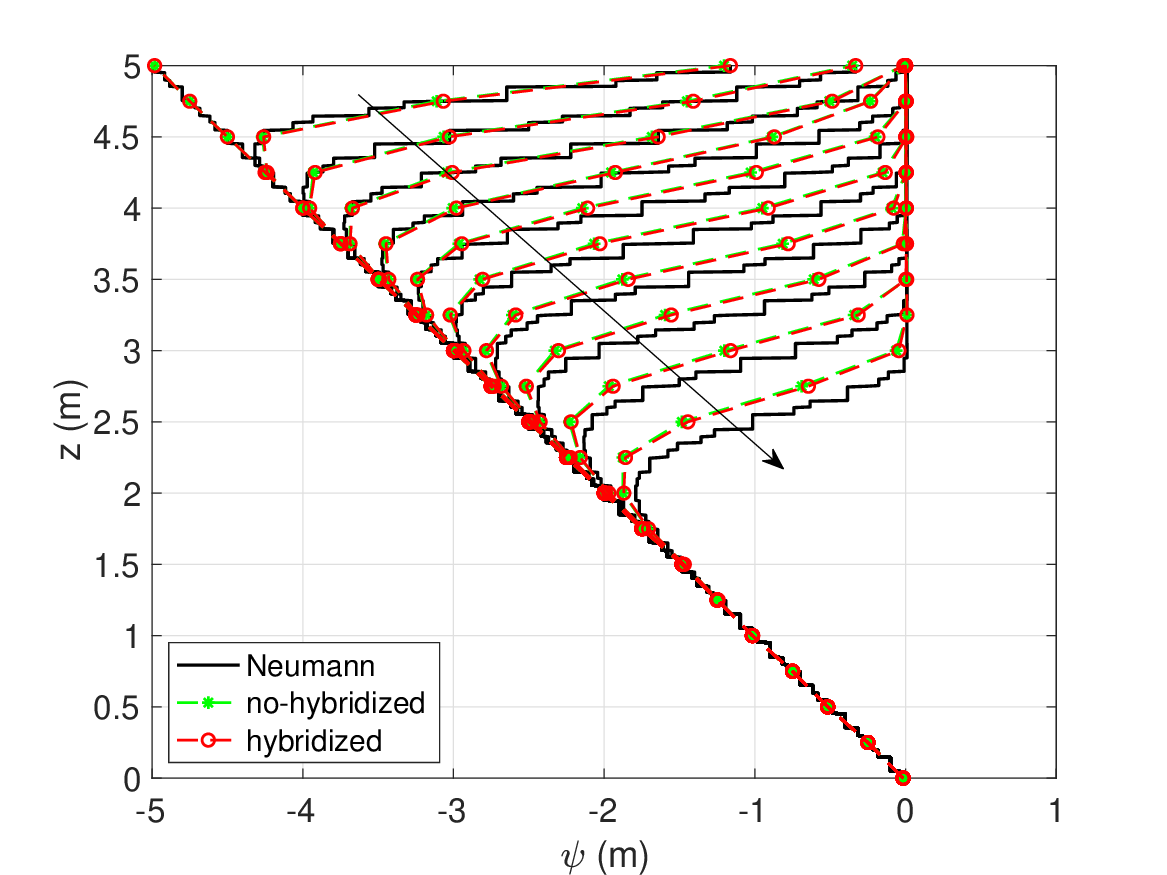})\hspace*{-.4cm}
\includegraphics[width=0.5\textwidth]{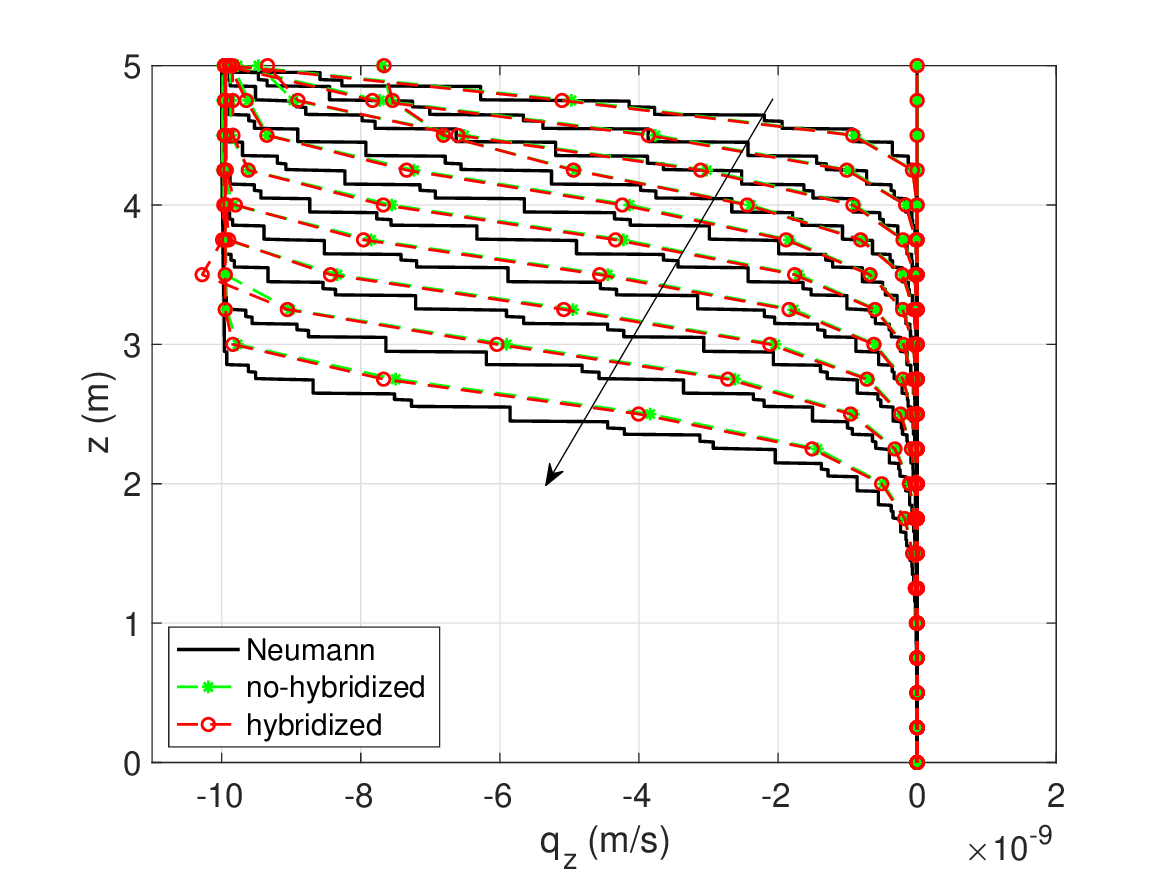}\vspace*{-.1cm}
\includegraphics[width=0.5\textwidth]{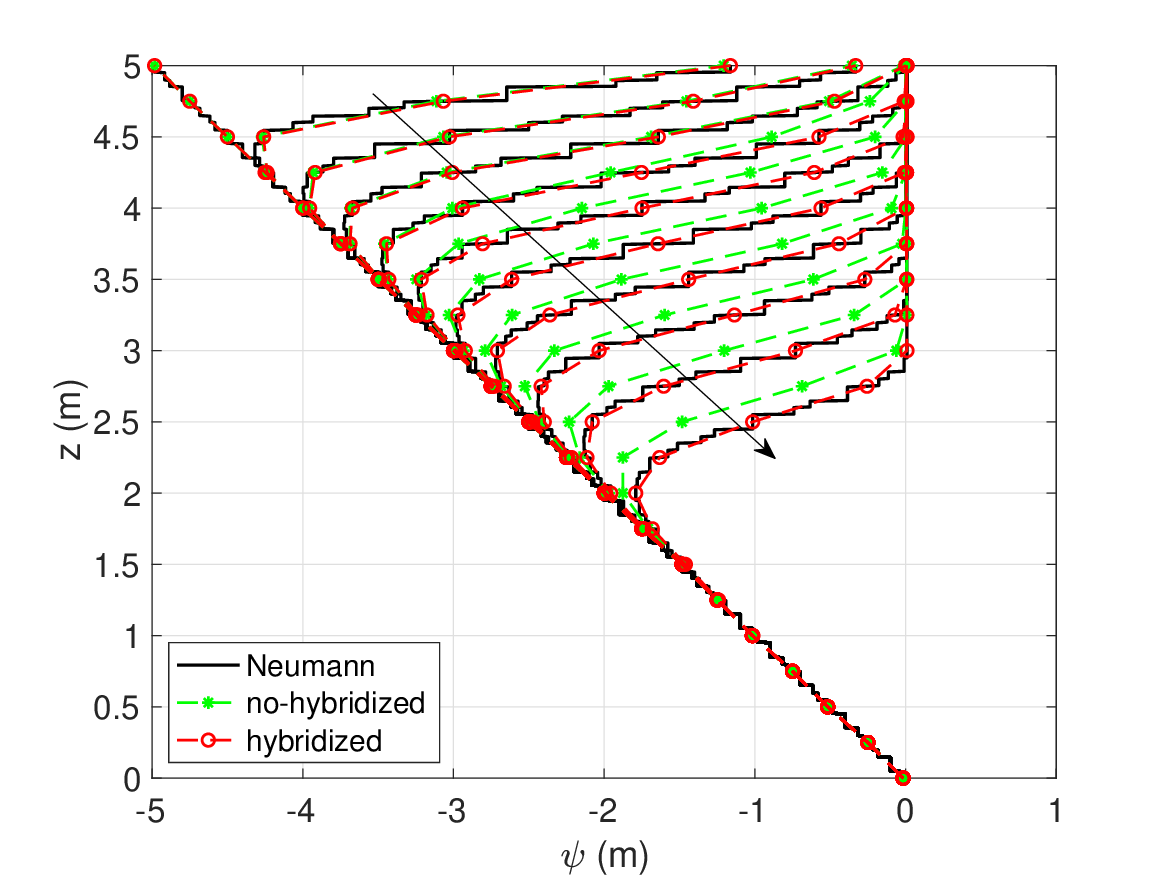}\hspace*{-.4cm}
\includegraphics[width=0.5\textwidth]{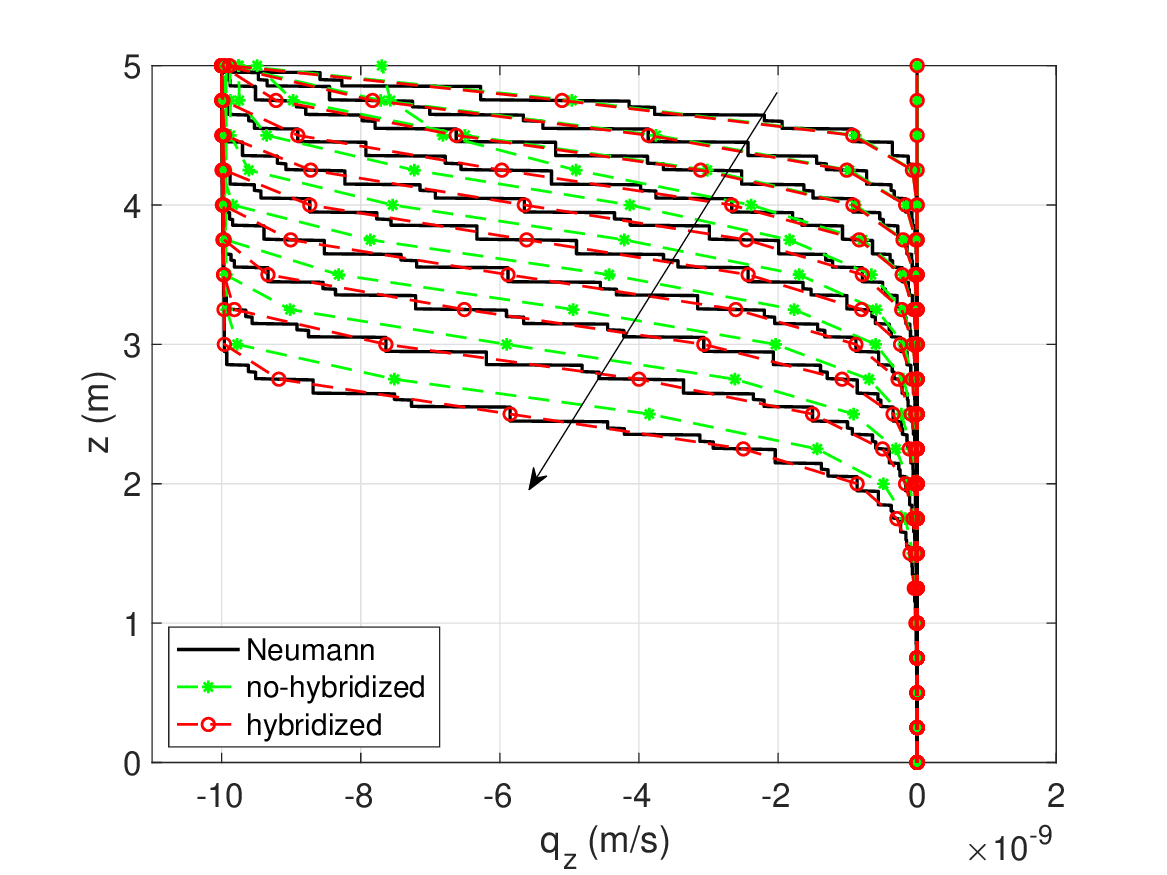}
\caption{Rectangular domain. Relaxation of the seepage conditions for certain soil types. Top row: results obtained by imposing $\epsilon=0$m on the two proposed discretizations, i.e., the no-hybridized and hybridized one. Bottom row: results obtained in case $\epsilon=10^{-2}$m. The curves are at times instants $0$h, $100$h, $200$h, \dots, $900$h, $1\,000$h.}
\label{fig:7}
\end{figure}
This relaxation approach is similar to what has been done in fluid-structure interaction works~\cite{burman2020nitsche, formaggia2021xfem}. More specifically, we increase the largest pressure head allowed on the surface,  
replacing \eqref{eq:model1:3} with
\begin{equation}\label{eq:relax}
\psi_h^{n+1} \le \epsilon, \quad
Q_h^{n+1} \le 0, \quad
Q_h^{n+1}(\psi_h^{n+1} - \epsilon) = 0,
\end{equation}
where $\epsilon$ is a positive threshold beyond which the existence of significant water ponding is acknowledged. An example case where relaxing the seepage conditions improves the
solution is the following. We set the limit case $p/K_S=1$ so that we can compare the results with the exact solution given by the solely imposition of the Neumann boundary conditions on $\Gt$. We discretize the domain with a mesh spacing equal to $0.1$m, and set a final time $t_{\text{fin}}=1\,000$h with a time discretization $\Delta t = 1$h. The maximum number of nonlinear iteration is set equal to $20$. To solve internal nonlinearities in the Richards' equation we employ the $L$-scheme.
In Figure~\ref{fig:7} we report the results of our analysis. In the top row we show the results obtained in case $\epsilon=0$m for the three solutions considered (Neumann boundary condition, seepage conditions by non-hybridized method, seepage conditions by hybridized method). As one can notice, both the 
hybridized and non-hybridized discretizations present a delay with respect to the limit case obtained by imposing Neumann boundary conditions on $\Gt$. 
In the bottom row we show the results obtained by setting $\epsilon=10^{-2}$m. As the reader can notice, while the solution obtained by the non-hybridized scheme does not change significantly, the hybridized version completely compensates the delay accumulated with respect to the Neumann solution.

\subsubsection{Sensitivity analysis on the choice of the soil parameters}\label{sec:sens:soil_type}
The goal of this test is to show that our approaches can deal with various kind of soils even in the most relevant case for the application, i.e., $p/K_S>1$ (so far, we have discussed this case only for clayey soils). Thus, we fix $p/K_S=10$ and we consider either silty or sandy soil (note the test in Section \ref{sec:rel:soil_type} was considering silt but in the case $p/K_S=1$).

The key observation of this test is that we found that the $L$-scheme converges to a 
non-physical solution for the sandy soil (not reported here), such that using the Newton method is necessary for this kind of soil. Conversely, both the $L$-scheme and the Newton method were found to be adequate to compute the solution in the case of silty soil (within the constraint of the
time-step for Newton discussed in Section \ref{sec:sec_non-linearities}). The results reported here for the silty soil are obtained with the $L$-scheme.

\begin{figure}[h!]
\centering
\includegraphics[width=0.5\textwidth]{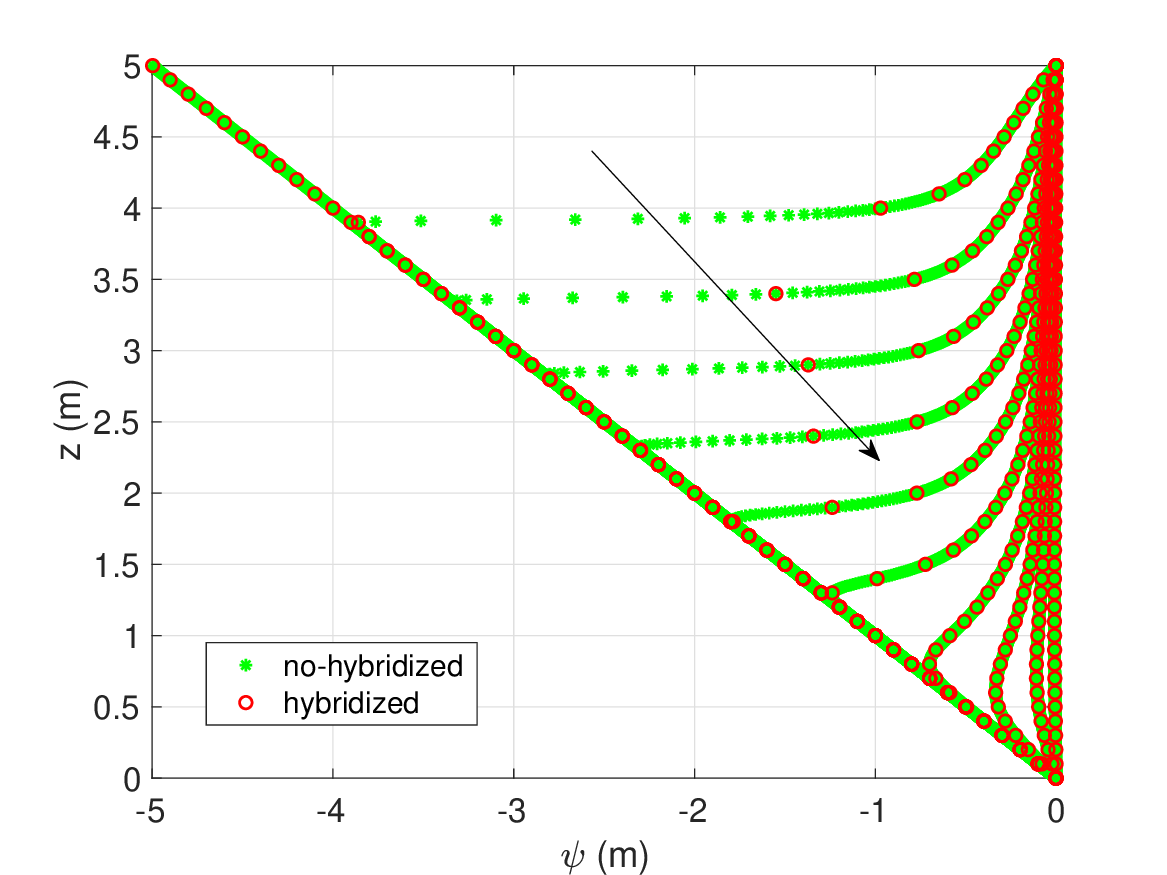}\hspace*{-.4cm}
\includegraphics[width=0.5\textwidth]{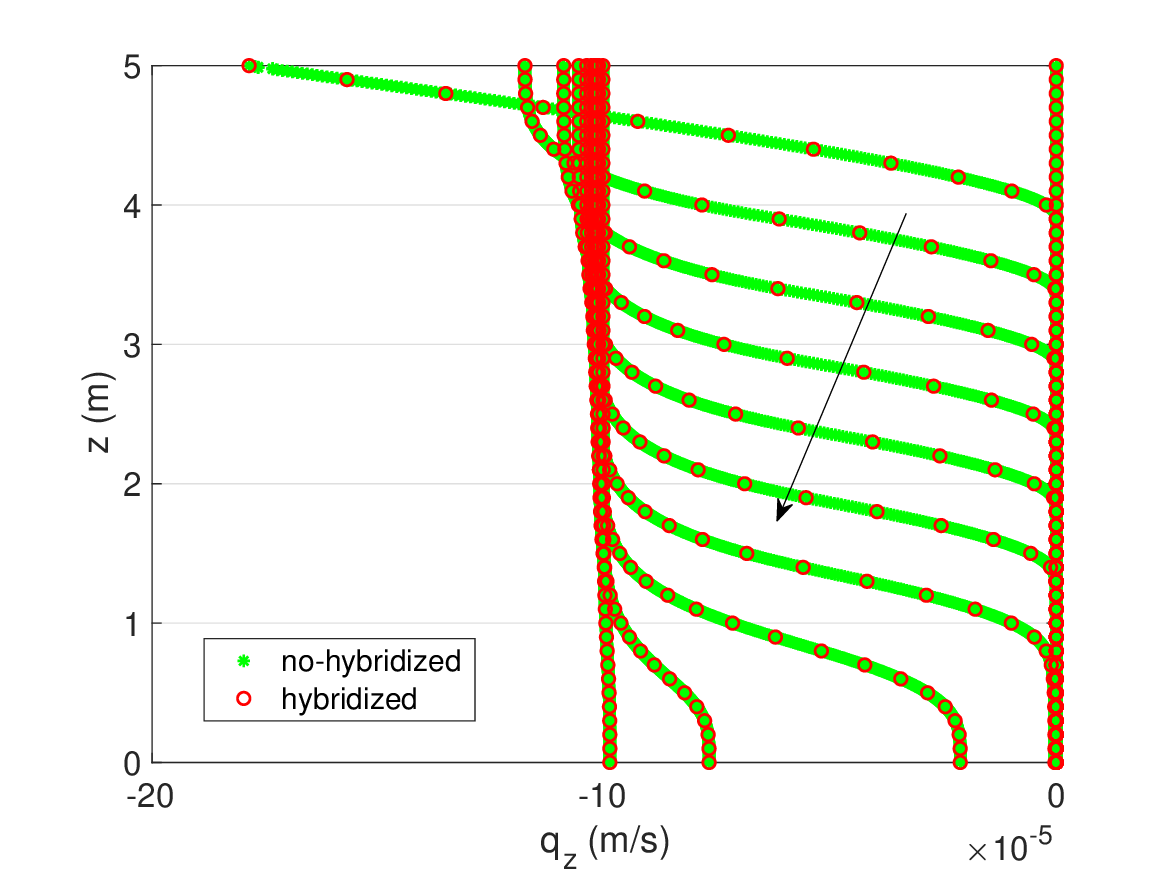}
\caption{Rectangular domain. Sensitivity analysis on the choice of the soil parameters. Sandy soil. The curves are at times instants $0$h, $0.5$h, $1$h, \dots, $5.5$h, $6$h. }
\label{fig:3}
\end{figure}

We set a final time $t_{\text{fin}}=6$h, and a time step $\Delta t = 1\,800$s.
The domain $\Omega$ is discretized with a mesh size of $0.01$m. 
For the case of sandy soil, we do not consider any relaxation on the seepage conditions, i.e., we set $\epsilon=0$ in Equation~\eqref{eq:relax}. 
In Figure~\ref{fig:3}, we report the results of this simulation, where the curves are temporally spaced by $1\,800$s. We particularly point out how both proposed schemes avoid ponding on $\Gt$, since the pressure on $\Gt$ is null and the flux is non-zero. The relative difference between the two computed solutions in $L^{\infty}(\Omega)$-norm is of the order of $10^{-15}$.

\begin{figure}[h!]
\centering
\includegraphics[width=0.5\textwidth]{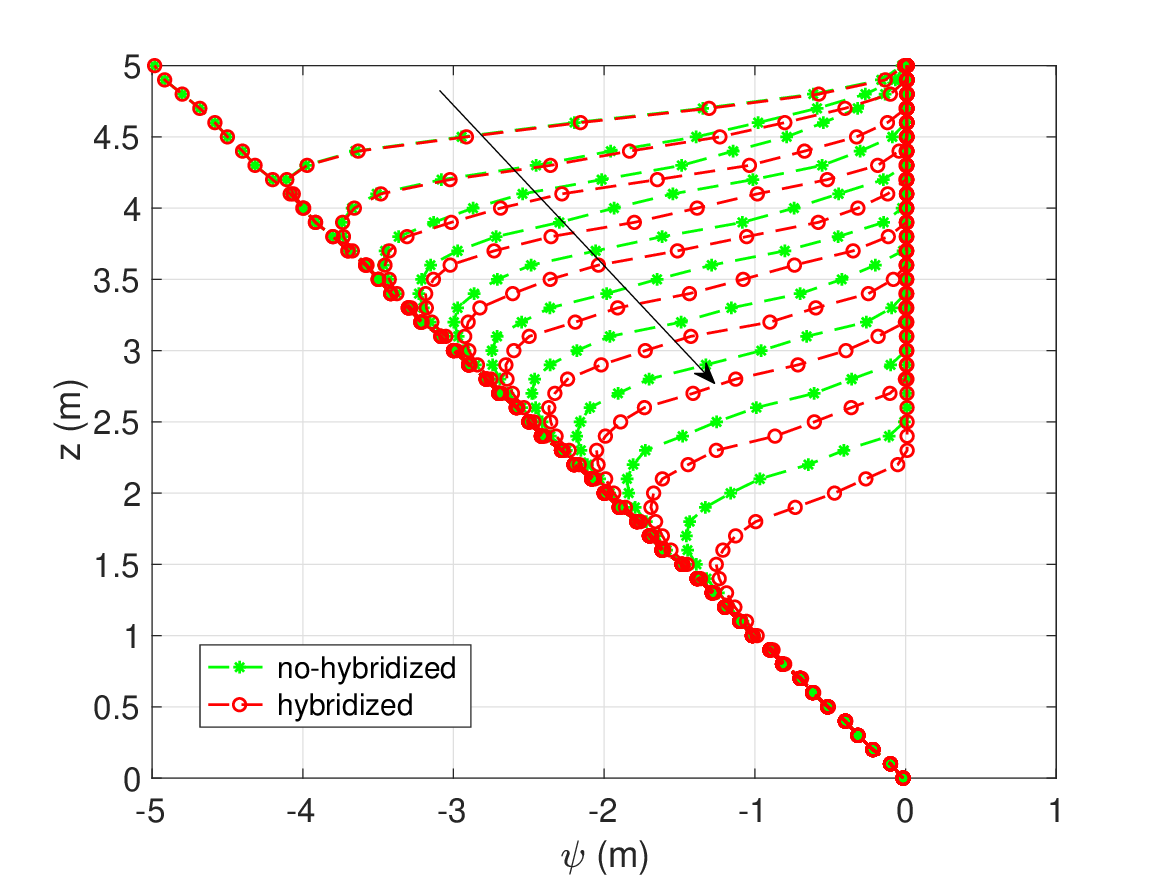}\hspace*{-.4cm}
\includegraphics[width=0.5\textwidth]{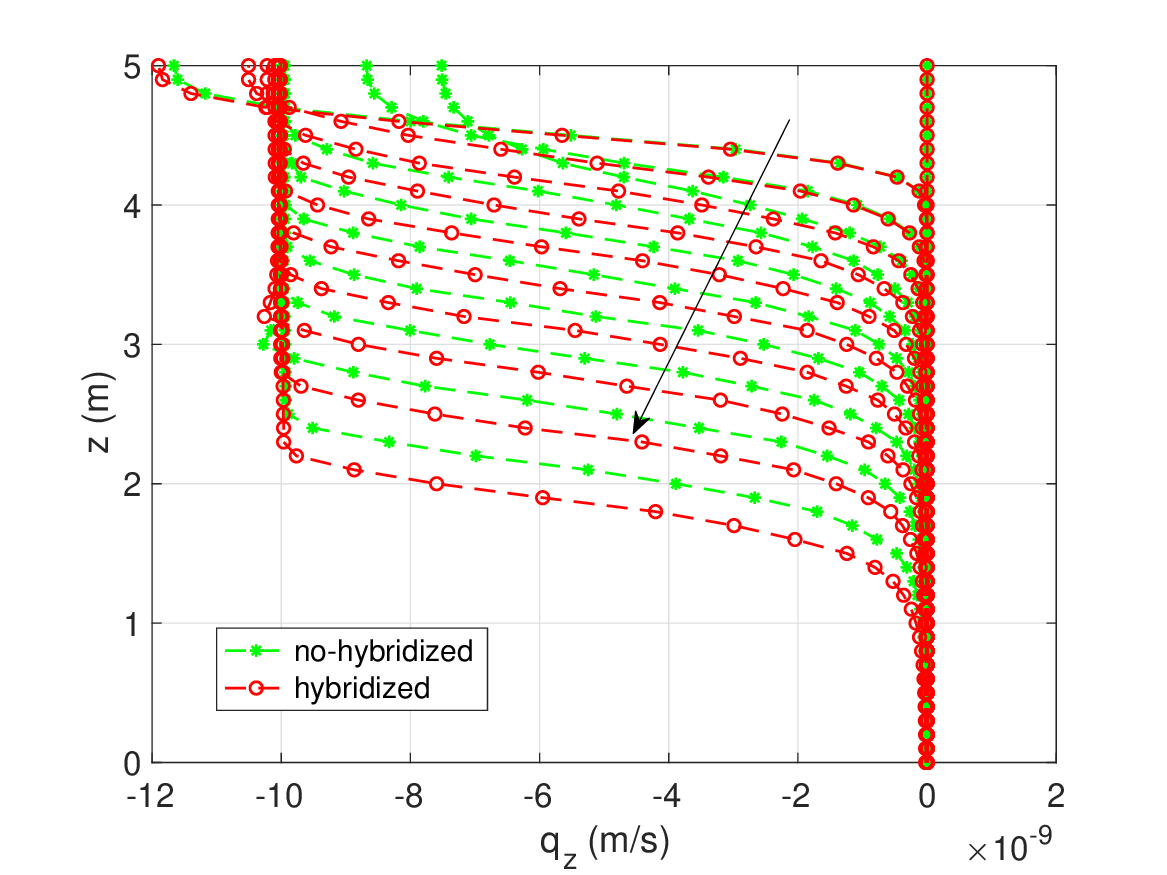}
\caption{Rectangular domain. Sensitivity analysis on the choice of the soil parameters. Silty soil. The curves are at times instants $0$h, $100$h, $200$h, \dots, $900$h, $1\,000$h.}
\label{fig:4}
\end{figure}

When considering the silty soil, we set the relaxation parameter $\epsilon$ in Equation \eqref{eq:relax} to $\epsilon=10^{-2}$m. Figure~\ref{fig:4} shows the results of this analysis. Again, we note that both proposed discretizations can handle the seepage conditions, as the pressure head never rises above $\epsilon$. However, it is evident that 
the temporal evolution of the non-hybridized scheme is slower than that of the 
hybridized scheme, as already observed in Section~\ref{sec:rel:soil_type}
for the case $p/K=1$.

\subsection{Artificial slope domain}\label{sec:ideal_slope}

\begin{figure}[h!]
\centering
\tikz{
\draw
    (0,0) node {\includegraphics[width=.7\textwidth]{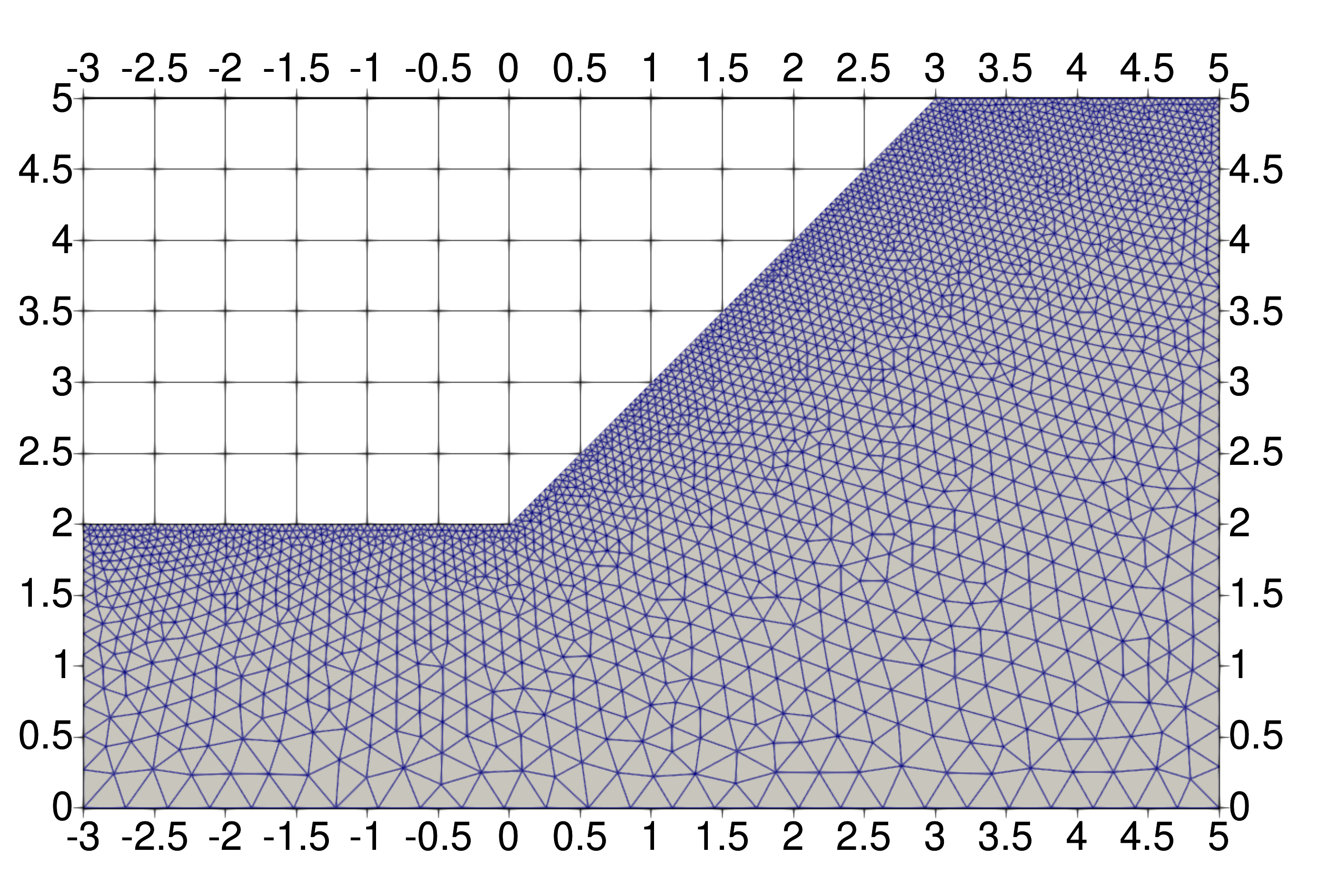}
};
\draw[<->,line width=0.9]
    (-.375\textwidth,-.23\textwidth)++(0.5,0.5) node[left,scale=1] {$z\,(m)$} |-++
         (0.5,-0.5) node[below,scale=1] {$x\,(m)$};
}
\caption{Artificial slope domain $\Omega$ and the computational mesh. }
\label{fig:8}
\end{figure}

In this section we face a so-called \emph{artificial slope} geometry. The aim of this test, along with the one presented in the next section, is to demonstrate the ability of the proposed discretizations to handle various domain complexities.
We consider the domain $\Omega$ depicted in Figure~\ref{fig:8}, which features a characteristic mesh size of $0.05$m on $\Gt$ and $0.3$m on $\Gb$. The total number of mesh elements is $5\,189$.

We consider a uniform and constant rainfall scenario with a ratio $p/K_S=10$ and a domain filled with clayey soil. We perform two tests with different initial conditions: a uniform initial condition corresponding to a pressure head of $\psi(x,z,0) = -20$m throughout the domain, and a hydrostatic initial condition with the water table being at toe of the slope, i.e., $\psi(x,z,0) = (2 - z)$m. In both cases, we impose a Dirichlet boundary condition on $\Gb$ to ensure compatibility with the initial data.
In particular, the second variant of this test is designed to demonstrate how the proposed method can automatically locate the position of the exit point on the potential seepage face $\Gt$. 
In the first test, the phreatic line is expected to advance from the top to the bottom of the domain, whereas in the second test, we anticipate the phreatic line to intersect the boundary $\Gt$ at some point during its evolution, thereby identifying the exit point. We perform the simulations up to $20$h with a time step of $5$h. We solve the Van Genuchten nonlinearities using the $L$-scheme with a tolerance for the stopping criteria of $\epsilon_A = 10^{-5}$. Below, we report only the results of the non-hybridized scheme with $\gamma_0 = 10^{-10}$, as we do not observe significant differences between the hybridized and non-hybridized solution schemes.
\begin{figure}[h!]
\centering
\tikz{
\draw
    (0,0) node {\includegraphics[width=.66\textwidth]{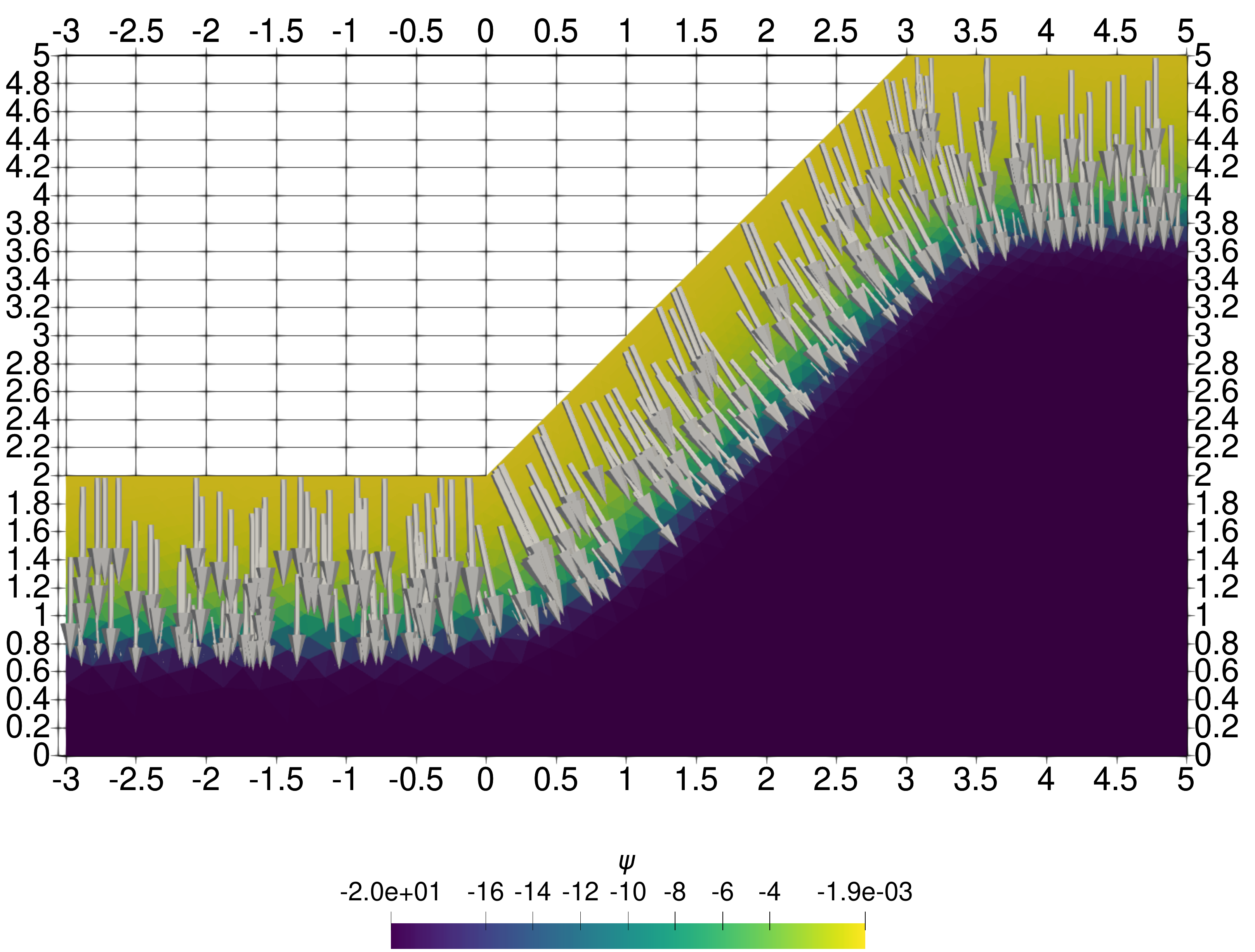}};
\draw[<->,line width=0.9]
    (-.375\textwidth,-.18\textwidth)++(0.5,0.5) node[left,scale=1] {$z\,(m)$} |-++
         (0.5,-0.5) node[below,scale=1] {$x\,(m)$};
}
\caption{Artificial slope domain. First test. We report the isolines of the pressure height $\psi$ together with the vector field $\vq$ at final time $20$h.  }
\label{fig:10}
\end{figure}

The first test we performed took roughly $23$s to complete, and required around $60$ iterations of the nonlinear solver to converge for each time step. In Figure~\ref{fig:10}, we show the results of the simulation. We present the isolines of the pressure head $\psi$ and the vectors of $\vq$ at the final time. As the reader can observe from the isolines, 
water infiltrates from the top of the domain, and a wetting front parallel to the ground surface develops, moving from the top to the bottom of the domain (look at the vectors representative of the vector flux $\vq$) \cite{sun1998analysis, lee2009simple, zhang2011stability}. We finally note that how the pressure head on the top boundary $\Gt$ approaches zero. 

\begin{figure}[h!]
\tikz{
\draw (0,0) node {
\begin{minipage}{\textwidth}
\centering
\includegraphics[width=.45\textwidth]{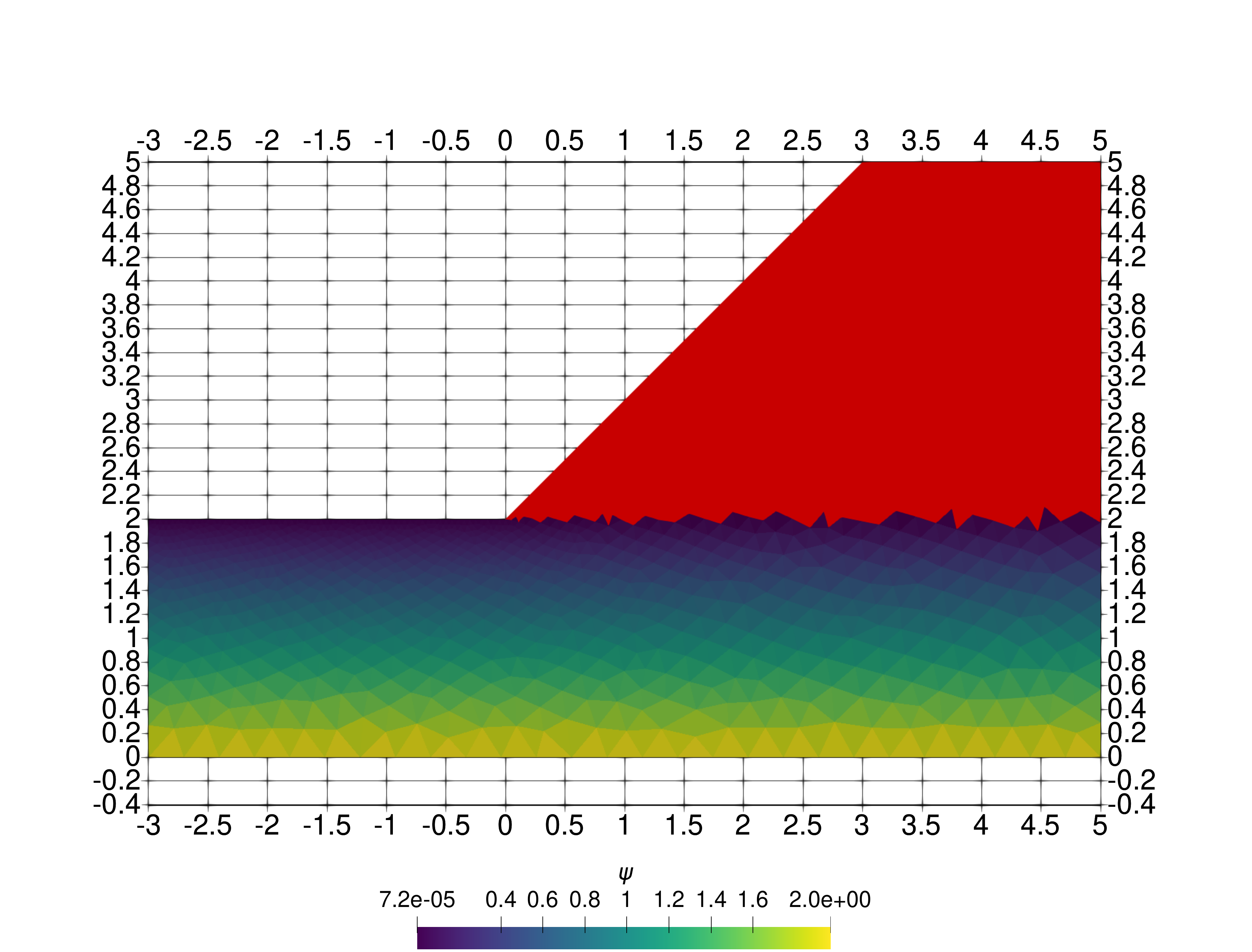}
\includegraphics[width=.45\textwidth]{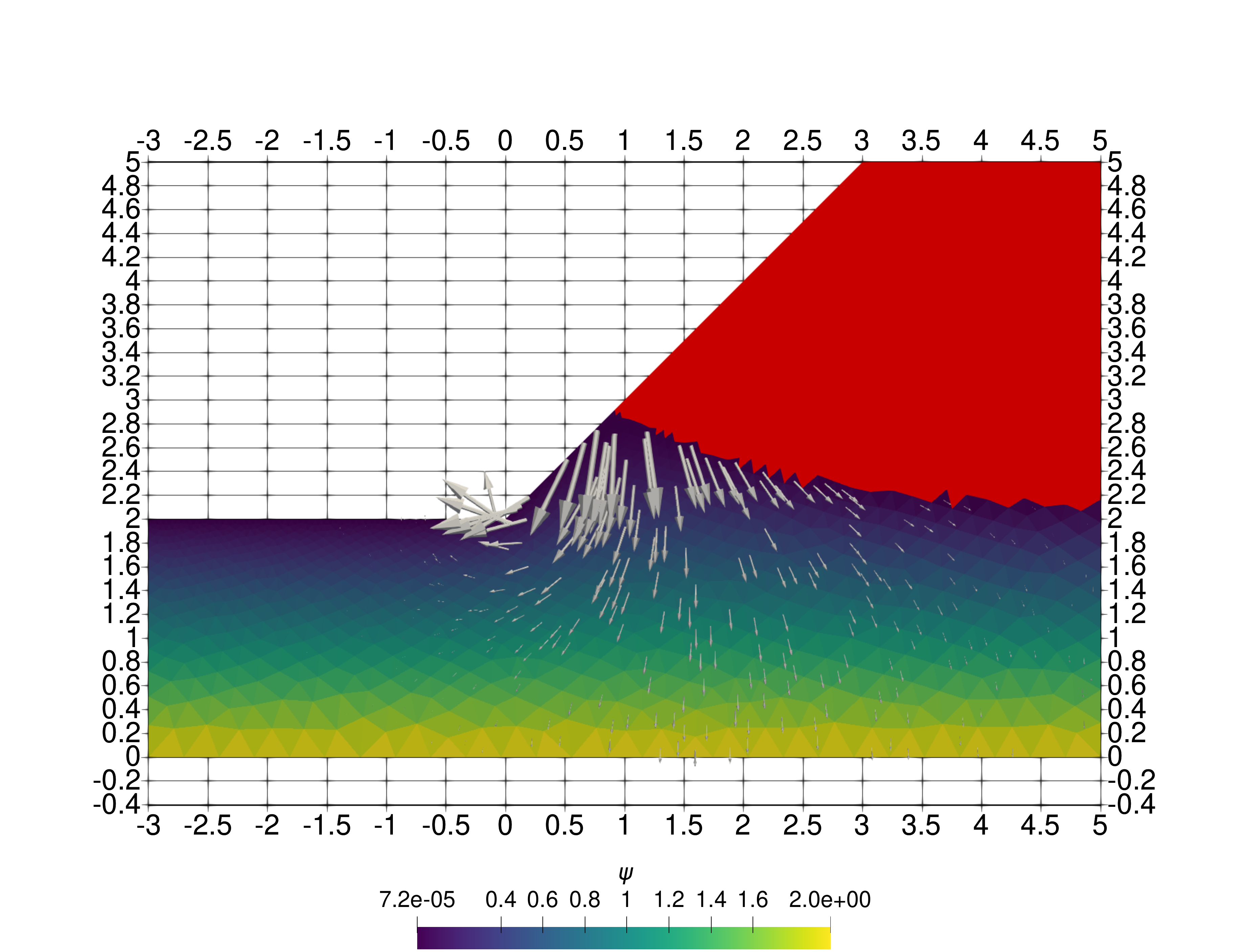}\\
\includegraphics[width=.45\textwidth]{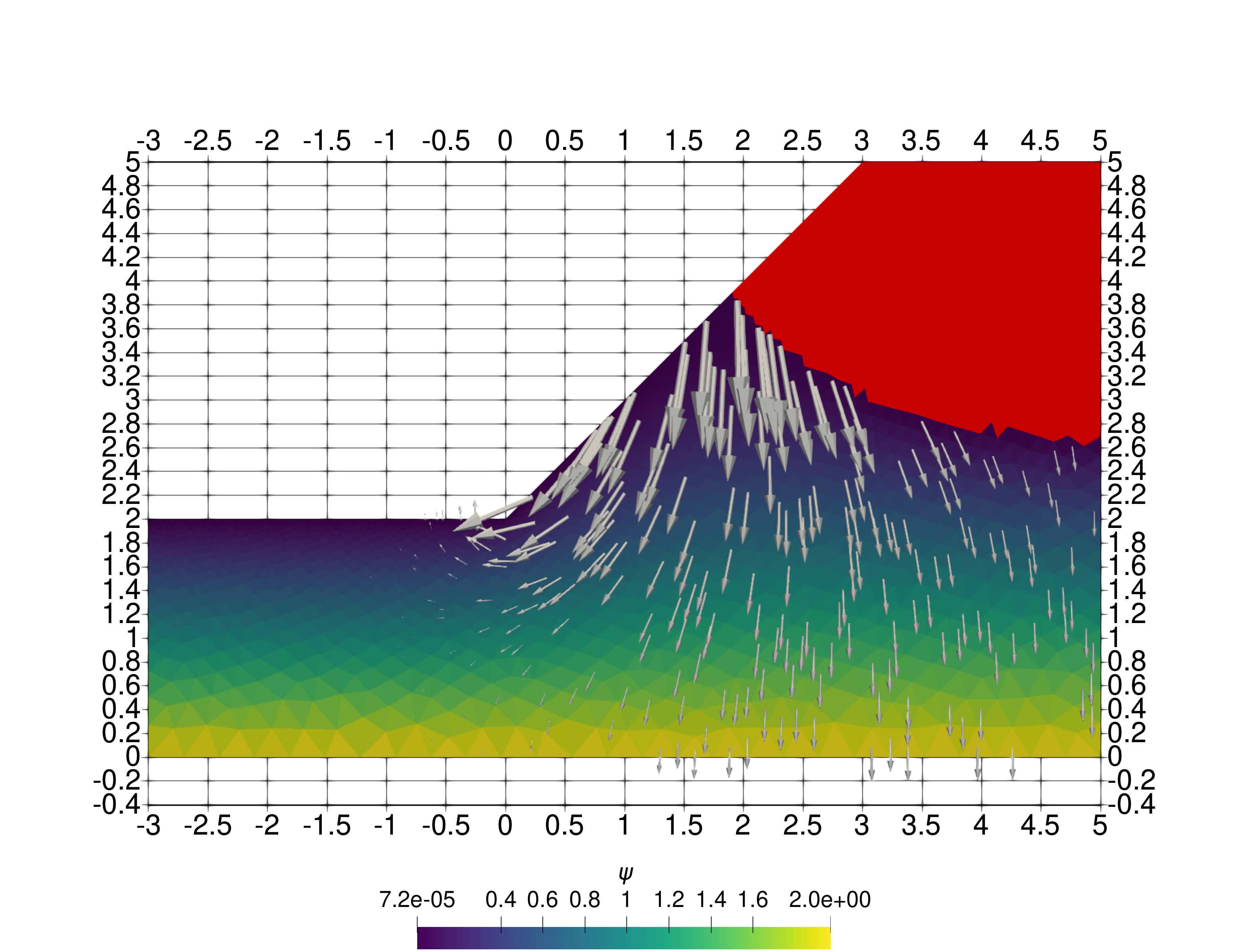}
\includegraphics[width=.45\textwidth]{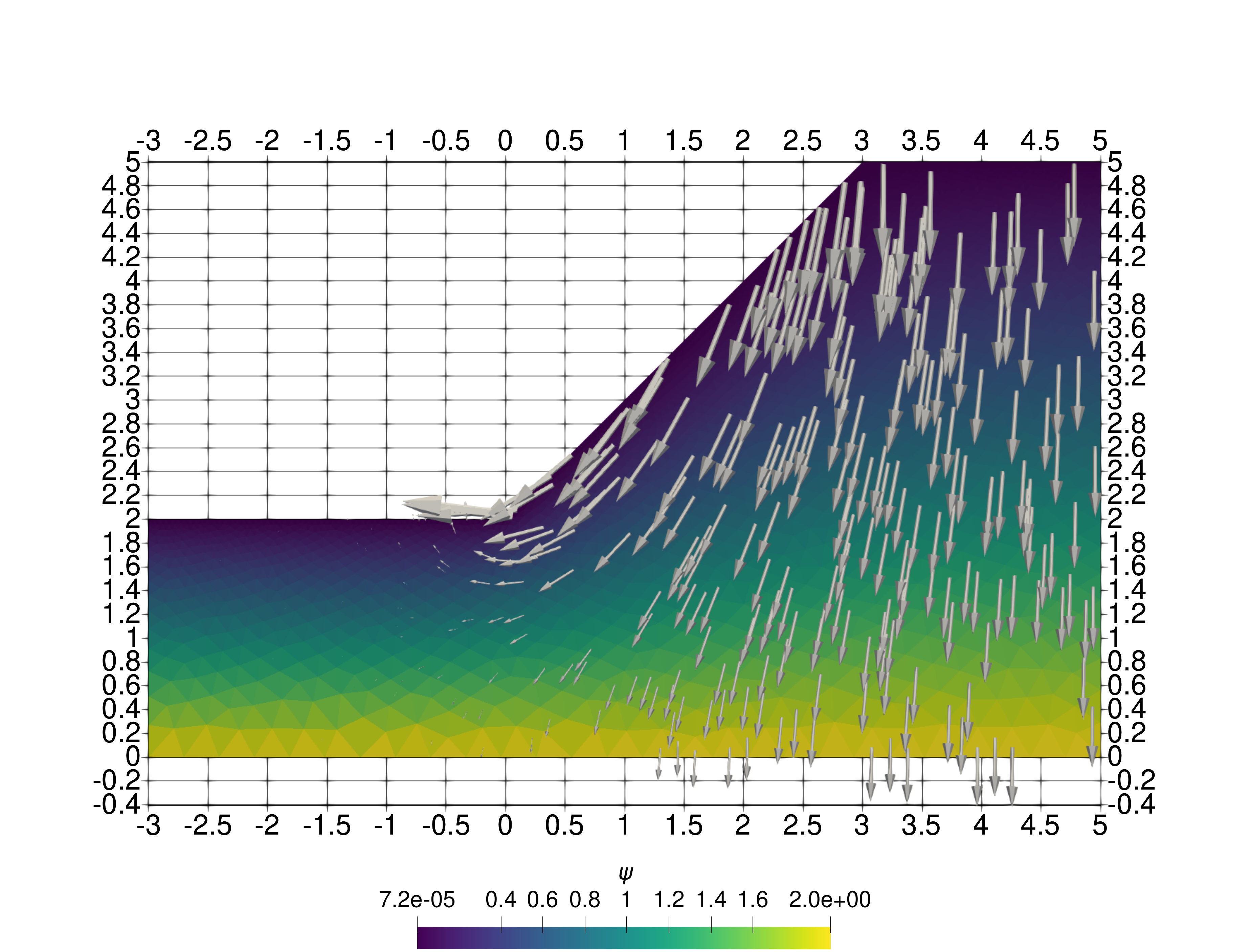}
\end{minipage}
};
\draw[<->,line width=0.9]
    (-.5\textwidth,-.325\textwidth)++(0.5,0.5) node[left,scale=1] {$z\,(m)$} |-++
         (0.5,-0.5) node[below,scale=1] {$x\,(m)$};}
\caption{Artificial slope domain. Second test. Top left panel is at initial time, top left panel is at time $5$h, bottom left is at time $10$h, bottom right panel is at time $20$h. We report the isolines of the pressure height $\psi$ together with the vector field $\vq$.  }
\label{fig:9}
\end{figure}

The second test we performed took roughly $6$s to complete, and required around $20$ iterations of the nonlinear solver to converge for each time step. In Figure~\ref{fig:9}, we show the results we obtained. The reader can observe the evolution of the phreatic line as we isolate the domain corresponding to a positive pressure head $\psi$ 
(the region in red corresponds to the unsaturated zone).
We display the evolution of the phreatic line over the simulation time. The top left panel shows the initial condition, top right panel the state at $5$h, bottom left panel at $10$h, and bottom right panel at $20$h. The portion of the phreatic line within the domain evolves, starting from the position $z=2$m, and moves upwards until it fully coincides with the boundary $\Gt$ at the end of the simulation.
Observe that at the toe of the slope water exits the domain $\vq\cdot\vn>0$, which is consistent with the conditions listed in Equation~\eqref{eq:model1:3} since $Q\le0$ allows water to exit as long as the pressure head $\psi$ is null (complementarity condition); this is here correctly simulated.

\subsection{Natural slope domain}\label{sec:real_slope_DTM}
In this test we consider a realistic topography profile coming from a DTM. 
The slope we consider is taken from Ville San Pietro, which is a small village in Liguria (Italy) situated at about $500$m a.s.l. on the right of the Impero River, refer to~\cite{bovolenta2020geotechnical, report} for a geotechnical description of the site. 
The depth of the rock substrate, assumed to be impermeable, was interpolated by ordinary kriging \cite{cressie1988spatial} using information collected from eight boreholes located mainly in the central portion where a landslide is present and from about thirty light penetrometric tests distributed over a larger area and some seismic profiles.

\begin{figure}[h!]
\centering
\tikz{
\draw[<->,line width=0.9]
    (4,3)++(0.5,0.5) node[left,scale=1] {$z\,(m)$} |-++
         (0.5,-0.5) node[below,scale=1] {$x\,(m)$};
\draw
    (2,0) node[right] {\includegraphics[width=.8\textwidth]{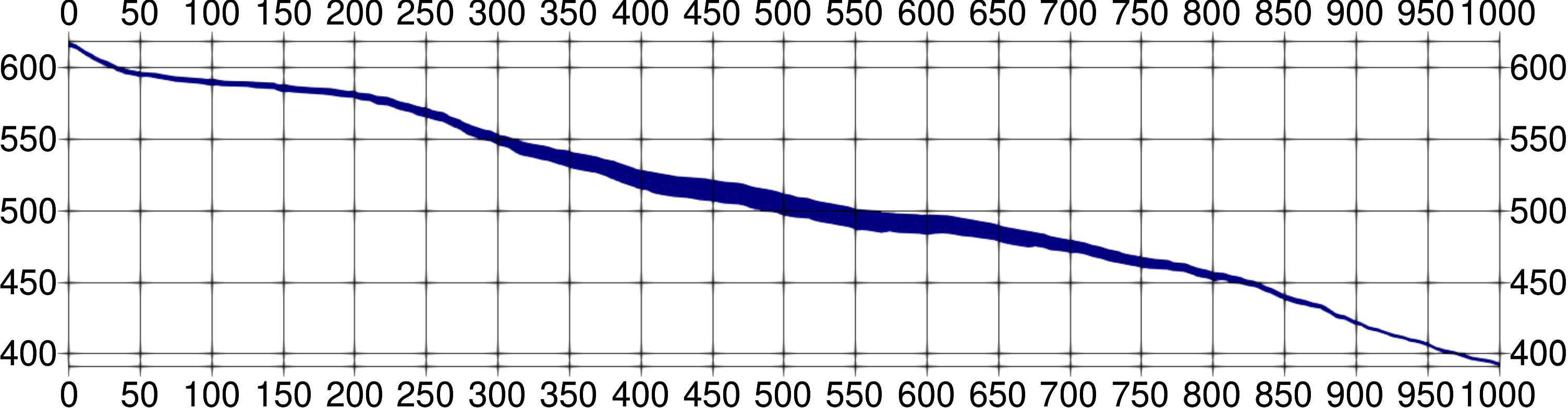}}
    (7,3) node[right]{\includegraphics[height=4cm]{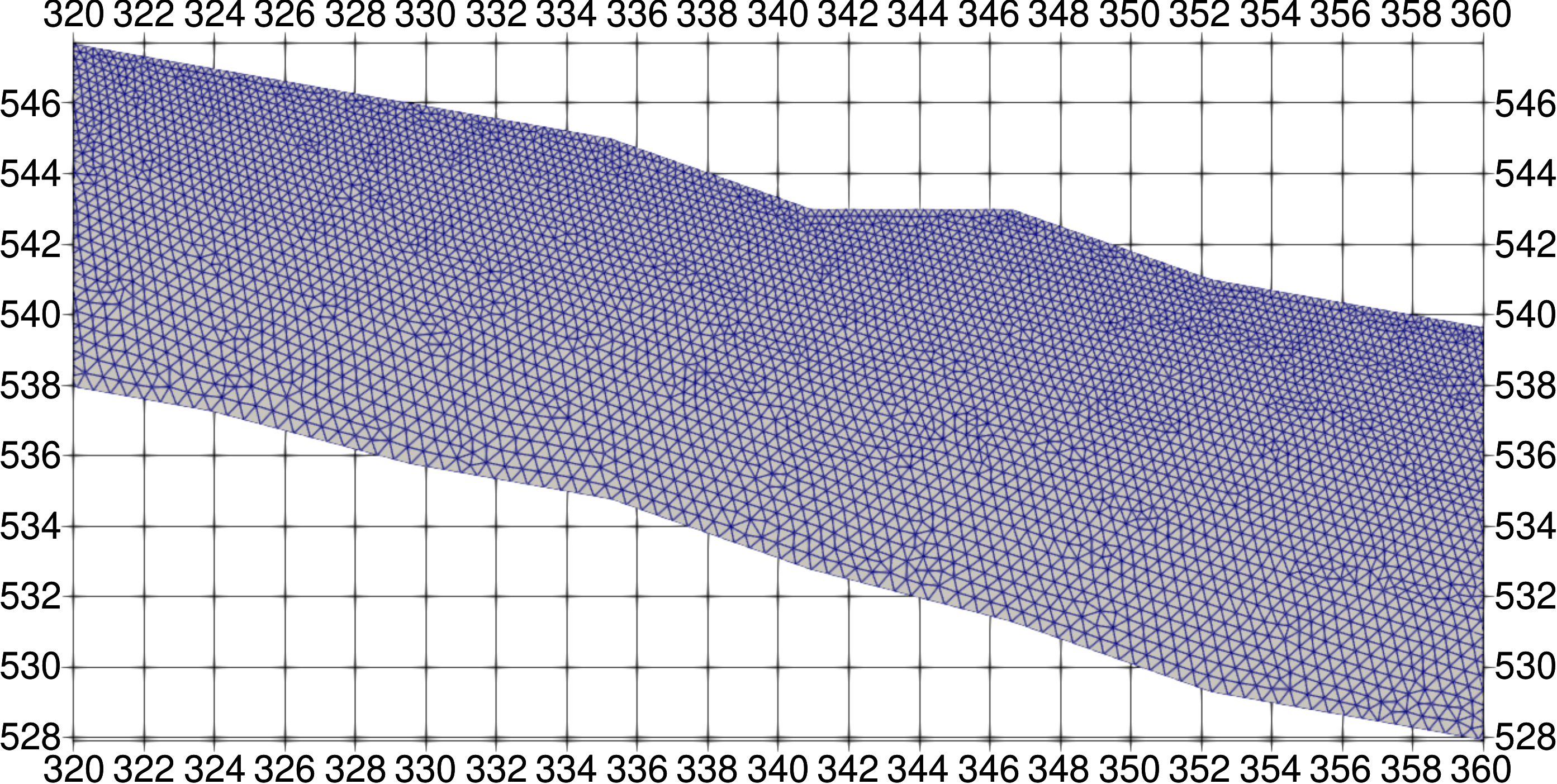}};
\fill[red,opacity=.5] (6.6015625,0.265625) rectangle (7.0625,0.5078125);
}
\caption{Natural slope domain $\Omega$ and the computational mesh. }
\label{fig:11}
\end{figure}

We present two scenarios, both characterized by the same soil and rainfall conditions ($p/K_S=10$) but with different initial and boundary condition on $\Gb$. While in the first test on $\Gb$ we enforce compatibility with initial data as boundary condition, on the second test we consider no-flux boundary condition. 
In both simulations, the Van~Genuchten nonlinearities are solved using a tolerance $\epsilon_A=10^{-5}$. In the first test we apply the $L$-scheme;
in the second test we use a combined method where $L$-scheme iterations are used until $\eta_{\text{lin}}^{k+1}< 10^{-3}$ and $k\ge100$, at which point we continue with Newton iterations to reach the $\epsilon_A=10^{-5}$ tolerance.

In Figure~\ref{fig:11} we show the domain $\Omega$ that entirely consists of clayey soil. 
The domain is discretized with a mesh size equal to $0.25$m at $\Gt$ and $0.5$m at $\Gb$.
The total number of triangular mesh elements is $173\,902$.
Below we present numerical results given by the hybridized approach. However, we point out that similar results are obtained by considering the non-hybridized version (not shown).

\begin{figure}[h!]
\centering
\tikz{

\draw (0.3,-1.5 ) node[above right] {\includegraphics[width=1\textwidth]{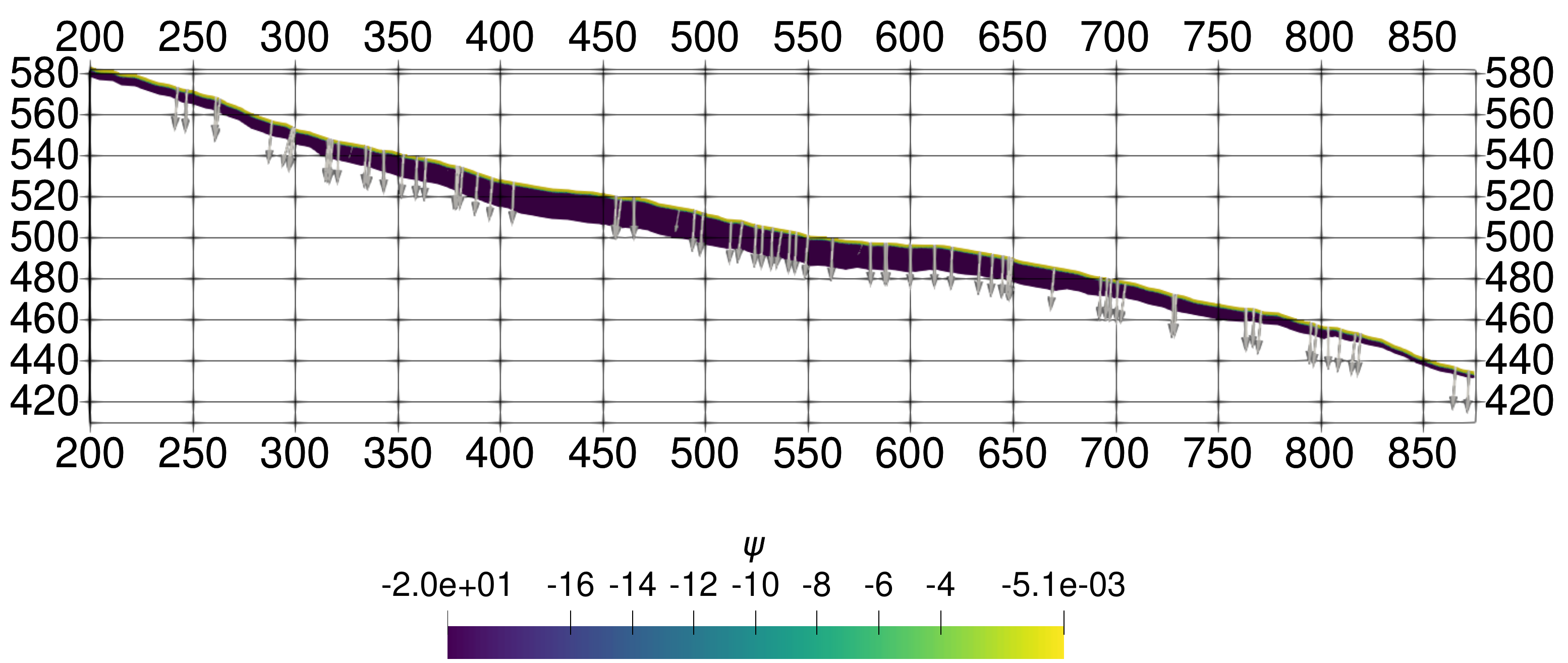}};
\draw[<->,line width=0.9]
    (0,0)++(0.5,0.5) node[left,scale=1] {$z\,(m)$} |-++
         (0.5,-0.5) node[below,scale=1] {$x\,(m)$};
}

\caption{Natural slope domain. First test case. Solution at final time $50$h. The arrows represent the vector flux $\vq$.  }
\label{fig:12}
\end{figure}

In a first test, we start from uniform initial conditions, i.e., $\psi(x,z,0)=-20$m throughout the entire domain. We set a final time of $50$h and a time step of $5$h. For each time step, the simulation took roughly $60$ iterations of the nonlinear solver to converge to the desired tolerance. The simulation took about $3\,800$ seconds to complete.
In Figure~\ref{fig:12}, we show the solution at $t=t_{\text{fin}}$. The results are in line with expectations. Indeed, the behavior is similar to what observed during the seepage analysis carried out in the rectangular domain and in the first test of Section~\ref{sec:ideal_slope}. We note the presence of an advancing wetting front starting from the boundary $\Gt$, which is corroborated by the arrows representing flux vectors $\vq$ pointing towards the boundary $\Gb$.

\begin{figure}[h!]
\centering
\tikz{
\draw (.3,-1.5) node[above right]{
\begin{minipage}{\textwidth}
\centering
\includegraphics[width=1\textwidth]{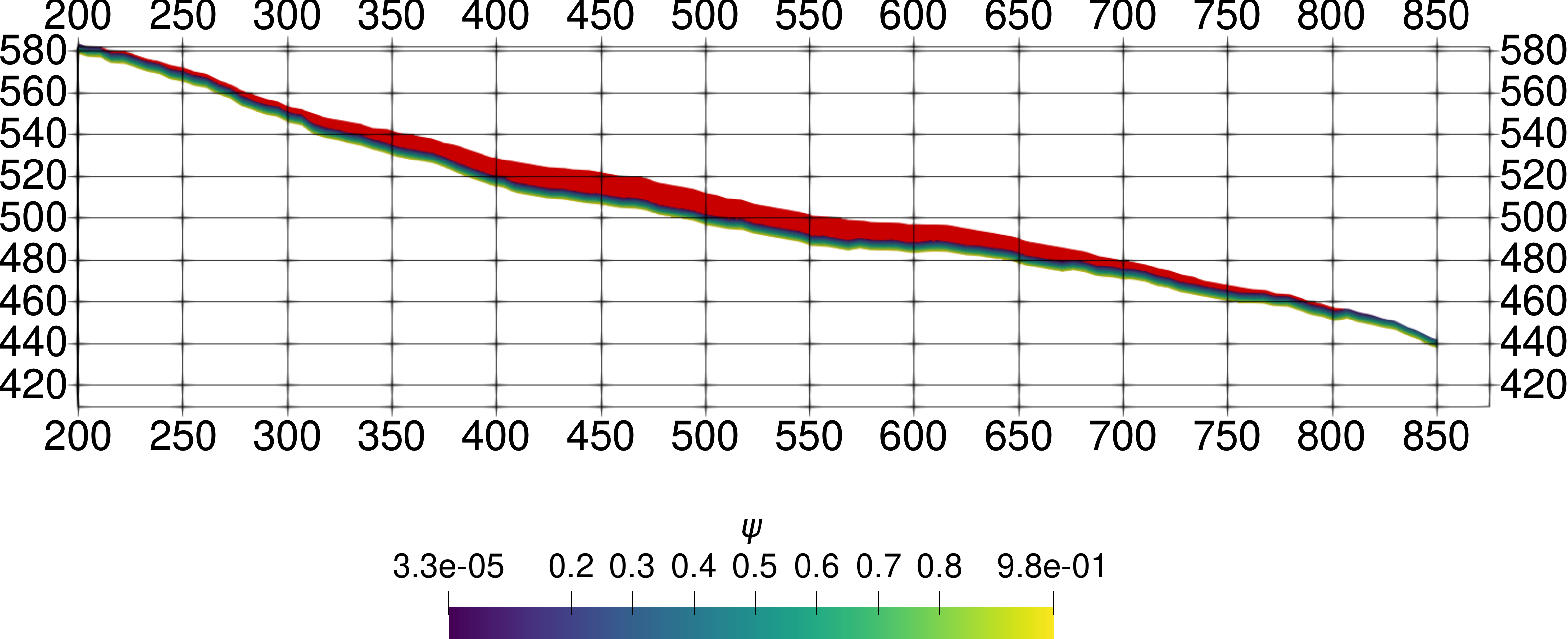}\\
\includegraphics[width=1\textwidth]{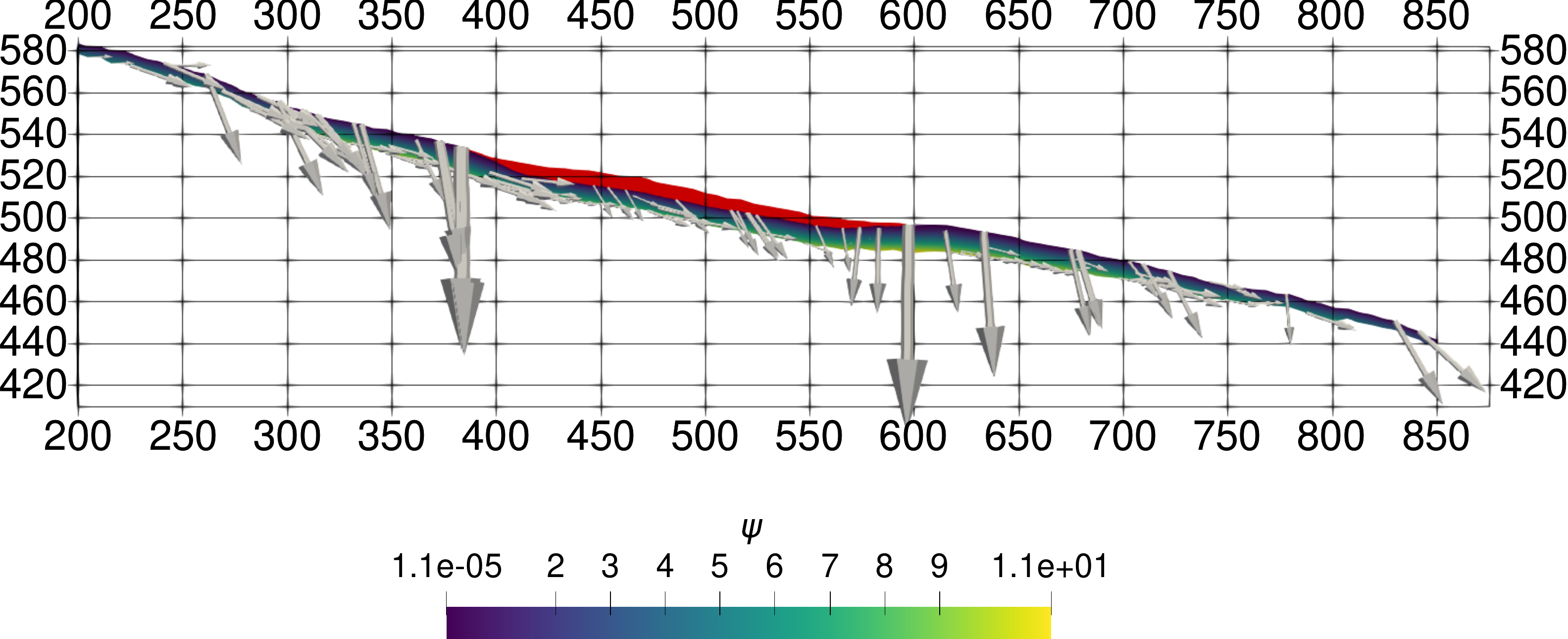}\\
\includegraphics[width=1\textwidth]{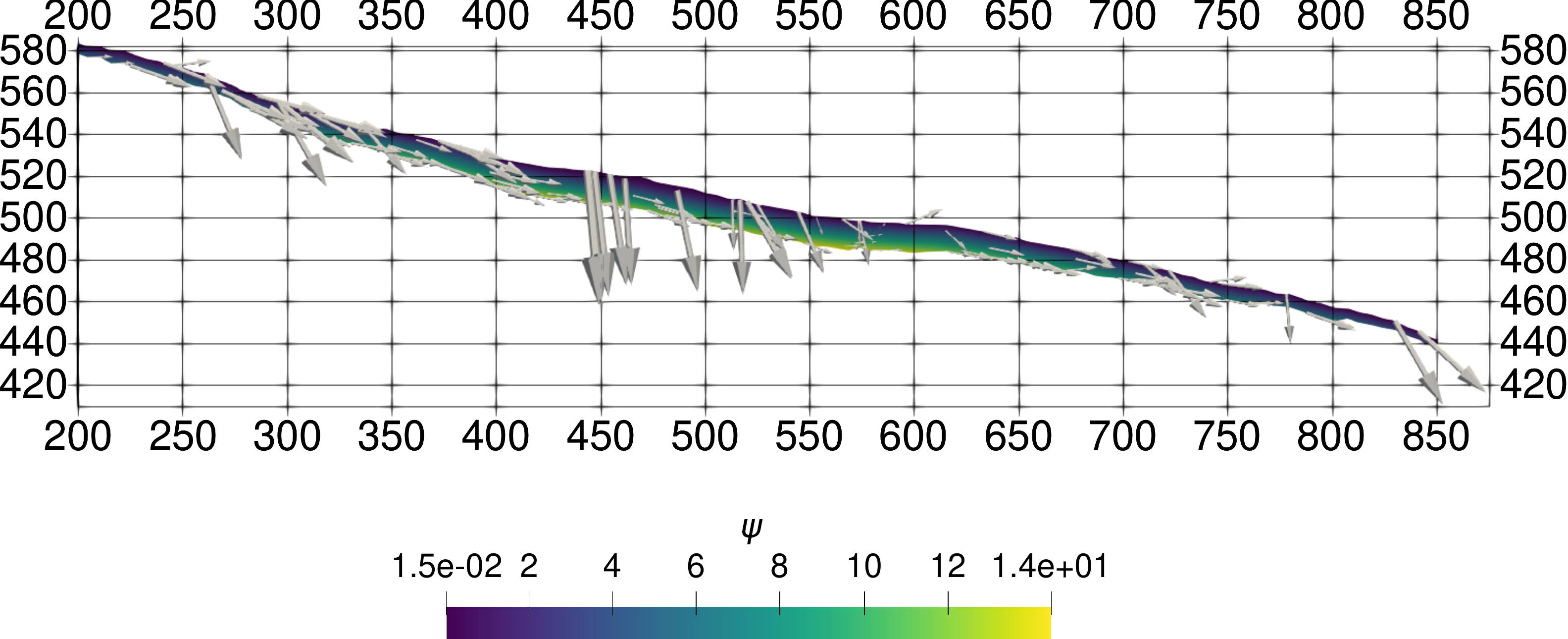}
\end{minipage}
};
\draw[<->,line width=0.9]
    (0,0)++(0.5,0.5) node[left,scale=1] {$z\,(m)$} |-++
         (0.5,-0.5) node[below,scale=1] {$x\,(m)$};}
\caption{Natural slope domain. Second test case. The red region corresponds to the fully unsaturated zone. The arrows represent the vector flux $\vq$. Top panel shows the initial condition. Middle and bottom panels represent the solution at times $10$h and $20$h, respectively. }
\label{fig:13}
\end{figure}

In a second test, we consider the initial condition $\psi(x,z,0)=1-\frac{2}{10}(z-z_{\text{bot}}(x))$, where $z_{\text{bot}}(x)$ is an equation identifying the boundary domain $\Gb$. This condition is chosen so that the phreatic line is inside the domain at $t=0$. As the simulation evolves, the wetting front advances in the soil and eventually reaches the impervious bedrock: from this moment on, 
a positive pressure builds up, leading to an increase of the water table level,
leading to the formation of seepage faces, that are correctly detected by the simulation. 
We set a final time of $40$h with a time step of $20$h, so that only two time steps are required.
The simulation took approximately $3\,043$ seconds to complete. The first time step required roughly $190$ nonlinear iterations to converge, while the second one about $150$. In Figure~\ref{fig:13}, we report the results of this analysis, presenting a zoomed-in view of a given spatial region.
Similarly to the test on the artificial slope (Section \ref{sec:ideal_slope}), the red region in the three panels represents the unsaturated zone, where $\psi\le 0$. 
The first panel shows the initial condition, while the middle and bottom panels show the solution at the two computed time steps. We particularly note how the proposed method can determine the location of all exit points on a real-case topography.
The presence of exit points on the slope can be deduced and partly confirmed by field surveys. There are several areas of vegetation typical of saturated soils and, in particular, for the exit points located approximately between $550$ and $600$ metres from the origin, there is a phreatimeter close to the profile. This has occasionally recorded the level of the water table coinciding with the topography.

\section{Conclusions}\label{sec:conclusion}
In this work we have proposed two discretization schemes that are able to impose the seepage boundary conditions for the Richards' equation, taking inspiration from the contact mechanics literature.
The first method relies on a strongly consistent penalization term, whereas the
second one is obtained by an hybridization approach, in which the value of the pressure on the
surface is treated as a separate set of unknown. Both methods 
can be understood as active-set methods, and 
are an improvement over existing heuristic algorithms that inevitably lead to water mass loss. Having a method that can automatically treat
both regimes of the seepage conditions is particularly advantageous when coupling the solver 
with other models for, e.g., meteorological data and slope stability.

We have discussed a sequence of numerical tests, that allowed us to showcase that both
approaches can effectively deal with the seepage conditions on computational domains of increasing
realism. In general, the hybridized method is preferable since the quality of the solution does not depend on the tuning of the penalization parameter, and possible delays in the temporal evolution
of the solution can be fixed by introducing suitable relaxations of the seepage condition.

\section*{Acknowledgements}

F.G., L.T., A.B., A.F. are members of the Gruppo Nazionale Calcolo Scientifico-Istituto Nazionale di Alta Matematica (GNCS-INdAM).

\section*{Declarations}

\subsection*{Declaration of competing interest}
The authors declare that they have no known competing financial interests or personal relationships that could have appeared to influence the work reported in this paper.

\subsection*{CRediT authorship contribution statement}
\textbf{Federico Gatti}: Conceptualization, Formal analysis, Investigation, Methodology, Software, Validation, Visualization, Writing – original draft, Writing – review \& editing. 
\textbf{Andrea Bressan}: Supervision, Writing – review \& editing, Validation. 
\textbf{Lorenzo Tamellini}: Supervision,  Writing – review \& editing, Validation, Funding acquisition. 
\textbf{Leonardo Maria Lalicata}: Supervision, Validation, Writing – review \& editing. 
\textbf{Simone Pittaluga}: Validation, Writing – review \& editing. 
\textbf{Domenico Gallipoli}: Supervision, Validation, Writing – review \& editing, Funding acquisition.
\textbf{Alessio Fumagalli}: Supervision, Writing – review \& editing, Validation, Writing – review \& editing.

\subsection*{Funding}
Funded by the European Union - NextGenerationEU and by the Ministry of University and Research (MUR), National Recovery and Resilience Plan (NRRP), Mission 4, Component 2, Investment 1.5, project “RAISE - Robotics and AI for Socio-economic Empowerment” (ECS00000035). F.G., L.T., A.B., S.P., D.G. are part of RAISE Innovation Ecosystem.

The present research is part of the activities of ``Dipartimento di Eccellenza 2023-2027'' of the Department of Mathematics of Politecnico di Milano, Italian Minister of University and Research (MUR), grant Dipartimento di Eccellenza 2023-2027.

\begin{figure}[h!]
\centering
\includegraphics[width=1\textwidth]{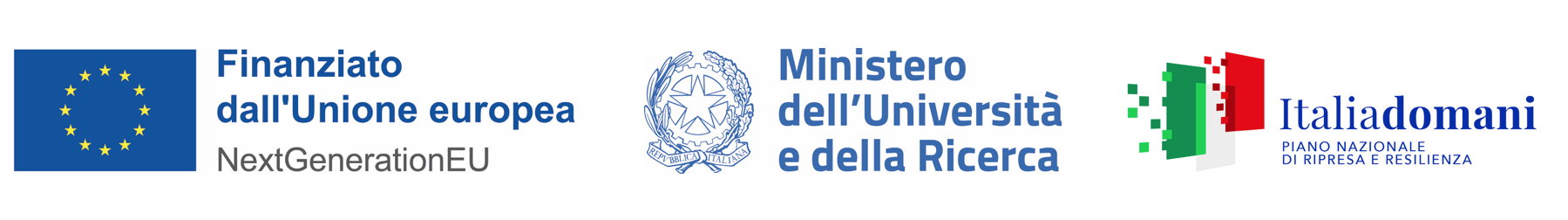}
\end{figure}


\bibliographystyle{elsarticle-num}
\bibliography{biblio.bib}

\begin{thebibliography}{10}
\expandafter\ifx\csname url\endcsname\relax
  \def\url#1{\texttt{#1}}\fi
\expandafter\ifx\csname urlprefix\endcsname\relax\def\urlprefix{URL }\fi
\expandafter\ifx\csname href\endcsname\relax
  \def\href#1#2{#2} \def\path#1{#1}\fi

\bibitem{rulon1985multiple}
J.~J. Rulon, R.~A. Freeze, Multiple seepage faces on layered slopes and their
  implications for slope-stability analysis, Canadian Geotechnical Journal
  22~(3) (1985) 347--356.

\bibitem{crosta1999slope}
G.~Crosta, C.~d. Prisco, On slope instability induced by seepage erosion,
  Canadian Geotechnical Journal 36~(6) (1999) 1056--1073.

\bibitem{orlandini2015evidence}
S.~Orlandini, G.~Moretti, J.~D. Albertson, Evidence of an emerging levee
  failure mechanism causing disastrous floods in {I}taly, Water Resources
  Research 51~(10) (2015) 7995--8011.

\bibitem{hirschfeld1973embankment}
R.~C. Hirschfeld, S.~J. Poulos, Embankment-dam engineering, 1973.

\bibitem{milligan2003some}
V.~Milligan, Some uncertainties in embankment dam engineering, Journal of
  Geotechnical and Geoenvironmental Engineering 129~(9) (2003) 785--797.

\bibitem{froude2018global}
M.~J. Froude, D.~N. Petley, Global fatal landslide occurrence from 2004 to
  2016, Natural Hazards and Earth System Sciences 18~(8) (2018) 2161--2181.

\bibitem{bishop1955use}
A.~W. Bishop, The use of the slip circle in the stability analysis of slopes,
  Geotechnique 5~(1) (1955) 7--17.

\bibitem{morgenstern1965analysis}
N.~U. Morgenstern, V.~E. Price, The analysis of the stability of general slip
  surfaces, Geotechnique 15~(1) (1965) 79--93.

\bibitem{pastor2021depth}
M.~Pastor, S.~M. Tayyebi, M.~M. Stickle, {\'A}.~Yag{\"u}e, M.~Molinos,
  P.~Navas, D.~Manzanal, A depth integrated, coupled, two-phase model for
  debris flow propagation, Acta Geotechnica (2021) 1--25.

\bibitem{he2023mpm}
K.~He, C.~Xi, B.~Liu, X.~Hu, G.~Luo, G.~Ma, R.~Zhou, {MPM-based mechanism and
  runout analysis of a compound reactivated landslide}, Computers and
  Geotechnics 159 (2023) 105455.

\bibitem{GATTI2024128525}
F.~Gatti, C.~{de Falco}, S.~Perotto, L.~Formaggia, A scalable well-balanced
  numerical scheme for the simulation of fast landslides with efficient time
  stepping, Applied Mathematics and Computation 468 (2024) 128525.

\bibitem{GATTI2024112798}
F.~Gatti, C.~{de Falco}, S.~Perotto, L.~Formaggia, M.~Pastor, A scalable
  well-balanced numerical scheme for the modeling of two-phase shallow granular
  landslide consolidation, Journal of Computational Physics 501 (2024) 112798.

\bibitem{freeze1971three}
R.~A. Freeze, Three-dimensional, transient, saturated-unsaturated flow in a
  groundwater basin, Water Resources Research 7~(2) (1971) 347--366.

\bibitem{rubin1968theoretical}
J.~Rubin, Theoretical analysis of two-dimensional, transient flow of water in
  unsaturated and partly unsaturated soils, Soil Science Society of America
  Journal 32~(5) (1968) 607--615.

\bibitem{neuman1975finite}
S.~P. Neuman, R.~A. Feddes, E.~Bresler, {Finite element analysis of
  two-dimensional flow in soils considering water uptake by roots: I. Theory},
  Soil Science Society of America Journal 39~(2) (1975) 224--230.

\bibitem{cooley1983some}
R.~L. Cooley, Some new procedures for numerical solution of variably saturated
  flow problems, Water Resources Research 19~(5) (1983) 1271--1285.

\bibitem{scudeler2017examination}
C.~Scudeler, C.~Paniconi, D.~Pasetto, M.~Putti, Examination of the seepage face
  boundary condition in subsurface and coupled surface/subsurface hydrological
  models, Water Resources Research 53~(3) (2017) 1799--1819.

\bibitem{chouly2013nitsche}
F.~Chouly, P.~Hild, {A Nitsche-based method for unilateral contact problems:
  numerical analysis}, SIAM Journal on Numerical Analysis 51~(2) (2013)
  1295--1307.

\bibitem{chouly2015symmetric}
F.~Chouly, P.~Hild, Y.~Renard, {Symmetric and non-symmetric variants of
  Nitsche’s method for contact problems in elasticity: theory and numerical
  experiments}, Mathematics of Computation 84~(293) (2015) 1089--1112.

\bibitem{chouly2017overview}
F.~Chouly, M.~Fabre, P.~Hild, R.~Mlika, J.~Pousin, Y.~Renard, {An overview of
  recent results on Nitsche’s method for contact problems}, in: Geometrically
  Unfitted Finite Element Methods and Applications: Proceedings of the UCL
  Workshop 2016, Springer, 2017, pp. 93--141.

\bibitem{chouly2018unbiased}
F.~Chouly, R.~Mlika, Y.~Renard, {An unbiased Nitsche’s approximation of the
  frictional contact between two elastic structures}, Numerische Mathematik 139
  (2018) 593--631.

\bibitem{raviart2006mixed}
P.-A. Raviart, J.-M. Thomas, A mixed finite element method for 2-nd order
  elliptic problems, in: Mathematical Aspects of Finite Element Methods:
  Proceedings of the Conference Held in Rome, December 10--12, 1975, Springer,
  2006, pp. 292--315.

\bibitem{brezzifortin1991}
F.~Brezzi, M.~Fortin, \href{https://doi.org/10.1007/978-1-4612-3172-1}{Mixed
  and hybrid finite element methods}, Vol.~15 of Springer Series in
  Computational Mathematics, Springer-Verlag, New York, 1991.
\newblock \href {https://doi.org/10.1007/978-1-4612-3172-1}
  {\path{doi:10.1007/978-1-4612-3172-1}}.
\newline\urlprefix\url{https://doi.org/10.1007/978-1-4612-3172-1}

\bibitem{tur2009mortar}
M.~Tur, F.~Fuenmayor, P.~Wriggers, A mortar-based frictional contact
  formulation for large deformations using {L}agrange multipliers, Computer
  Methods in Applied Mechanics and Engineering 198~(37-40) (2009) 2860--2873.

\bibitem{brunssen2007contact}
S.~Brun{\ss}en, S.~H{\"u}eber, B.~Wohlmuth, {Contact dynamics with Lagrange
  multipliers}, in: IUTAM Symposium on Computational Methods in Contact
  Mechanics: Proceedings of the IUTAM Symposium held in Hannover, Germany,
  November 5--8, 2006, Springer, 2007, pp. 17--32.

\bibitem{papadopoulos1998lagrange}
P.~Papadopoulos, J.~Solberg, A {L}agrange multiplier method for the finite
  element solution of frictionless contact problems, Mathematical and computer
  modelling 28~(4-8) (1998) 373--384.

\bibitem{burman2020nitsche}
E.~Burman, M.~A. Fern{\'a}ndez, S.~Frei, {A Nitsche-based formulation for
  fluid-structure interactions with contact}, ESAIM: Mathematical Modelling and
  Numerical Analysis 54~(2) (2020) 531--564.

\bibitem{burman2007stabilized}
E.~Burman, M.~A. Fern{\'a}ndez, {Stabilized explicit coupling for
  fluid--structure interaction using Nitsche's method}, Comptes Rendus.
  Math{\'e}matique 345~(8) (2007) 467--472.

\bibitem{burman2014unfitted}
E.~Burman, M.~A. Fern{\'a}ndez, {An unfitted Nitsche method for incompressible
  fluid--structure interaction using overlapping meshes}, Computer Methods in
  Applied Mechanics and Engineering 279 (2014) 497--514.

\bibitem{formaggia2021xfem}
L.~Formaggia, F.~Gatti, S.~Zonca, {An XFEM/DG approach for fluid-structure
  interaction problems with contact}, Applications of Mathematics 66~(2) (2021)
  183--211.

\bibitem{berge2020finite}
R.~L. Berge, I.~Berre, E.~Keilegavlen, J.~M. Nordbotten, B.~Wohlmuth, Finite
  volume discretization for poroelastic media with fractures modeled by contact
  mechanics, International Journal for Numerical Methods in Engineering 121~(4)
  (2020) 644--663.

\bibitem{cai1997control}
Z.~Cai, J.~Jones, S.~McCormick, T.~Russell, Control-volume mixed finite element
  methods, Computational Geosciences 1 (1997) 289--315.

\bibitem{aavatsmark2002introduction}
I.~Aavatsmark, An introduction to multipoint flux approximations for
  quadrilateral grids, Computational Geosciences 6 (2002) 405--432.

\bibitem{bianchi2022analysis}
D.~Bianchi, D.~Gallipoli, R.~Bovolenta, M.~Leoni, Analysis of unsaturated
  seepage in infinite slopes by means of horizontal ground infiltration models,
  G{\'e}otechnique (2022) 1--9.

\bibitem{van1980closed}
M.~T. Van~Genuchten, A closed-form equation for predicting the hydraulic
  conductivity of unsaturated soils, Soil science society of America journal
  44~(5) (1980) 892--898.

\bibitem{kikuchi1988contact}
N.~Kikuchi, J.~T. Oden, Contact problems in elasticity: a study of variational
  inequalities and finite element methods, SIAM, 1988.

\bibitem{wriggers2004computational}
P.~Wriggers, G.~Zavarise, {Computational contact mechanics. Encyclopedia of
  computational mechanics}, Solids and structures 2 (2004) 195--226.

\bibitem{ALART1991353}
P.~Alart, A.~Curnier, A mixed formulation for frictional contact problems prone
  to {N}ewton like solution methods, Computer methods in applied mechanics and
  engineering 92~(3) (1991) 353--375.

\bibitem{slodicka2002robust}
M.~Slodicka, A robust and efficient linearization scheme for doubly nonlinear
  and degenerate parabolic problems arising in flow in porous media, SIAM
  journal on scientific computing 23~(5) (2002) 1593--1614.

\bibitem{pop2004mixed}
I.~S. Pop, F.~Radu, P.~Knabner, {Mixed finite elements for the Richards’
  equation: linearization procedure}, Journal of computational and applied
  mathematics 168~(1-2) (2004) 365--373.

\bibitem{list2016study}
F.~List, F.~A. Radu, {A study on iterative methods for solving Richards’
  equation}, Computational Geosciences 20 (2016) 341--353.

\bibitem{bergamaschi1999mixed}
L.~Bergamaschi, M.~Putti, {Mixed finite elements and Newton-type linearizations
  for the solution of Richards' equation}, International journal for numerical
  methods in engineering 45~(8) (1999) 1025--1046.

\bibitem{lehmann1998comparison}
F.~Lehmann, P.~Ackerer, Comparison of iterative methods for improved solutions
  of the fluid flow equation in partially saturated porous media, Transport in
  porous media 31 (1998) 275--292.

\bibitem{stokke2023adaptive}
J.~S. Stokke, K.~Mitra, E.~Storvik, J.~W. Both, F.~A. Radu, An adaptive
  solution strategy for {R}ichards' equation, Computers \& Mathematics with
  Applications 152 (2023) 155--167.

\bibitem{mitra2019modified}
K.~Mitra, I.~S. Pop, A modified {L}-scheme to solve nonlinear diffusion
  problems, Computers \& Mathematics with Applications 77~(6) (2019)
  1722--1738.

\bibitem{radu2006newton}
F.~A. Radu, I.~S. Pop, P.~Knabner, Newton—type methods for the mixed finite
  element discretization of some degenerate parabolic equations, in: Numerical
  Mathematics and Advanced Applications: Proceedings of ENUMATH 2005, the 6th
  European Conference on Numerical Mathematics and Advanced Applications
  Santiago de Compostela, Spain, July 2005, Springer, 2006, pp. 1192--1200.

\bibitem{hueber2008primal}
S.~H{\"u}eber, G.~Stadler, B.~I. Wohlmuth, {A primal-dual active set algorithm
  for three-dimensional contact problems with Coulomb friction}, SIAM Journal
  on scientific computing 30~(2) (2008) 572--596.

\bibitem{wohlmuth2011variationally}
B.~Wohlmuth, Variationally consistent discretization schemes and numerical
  algorithms for contact problems, Acta Numerica 20 (2011) 569--734.

\bibitem{PyGeoN}
A.~Fumagalli, W.~M. Boon, E.~Ballini, compgeo-mox/pygeon: Pygeon 0.4 (apr
  2024).
\newblock \href {https://doi.org/10.5281/zenodo.10974537}
  {\path{doi:10.5281/zenodo.10974537}}.

\bibitem{keilegavlen2021porepy}
E.~Keilegavlen, R.~Berge, A.~Fumagalli, M.~Starnoni, I.~Stefansson, J.~Varela,
  I.~Berre, {Porepy: An open-source software for simulation of multiphysics
  processes in fractured porous media}, Computational Geosciences 25 (2021)
  243--265.

\bibitem{geuzaine2009gmsh}
C.~Geuzaine, J.-F. Remacle, {Gmsh: A 3-D finite element mesh generator with
  built-in pre-and post-processing facilities}, International journal for
  numerical methods in engineering 79~(11) (2009) 1309--1331.

\bibitem{2020SciPy-NMeth}
P.~Virtanen, R.~Gommers, T.~E. Oliphant, M.~Haberland, T.~Reddy, D.~Cournapeau,
  E.~Burovski, P.~Peterson, W.~Weckesser, J.~Bright, S.~J. {van der Walt},
  M.~Brett, J.~Wilson, K.~J. Millman, N.~Mayorov, A.~R.~J. Nelson, E.~Jones,
  R.~Kern, E.~Larson, C.~J. Carey, {\.I}.~Polat, Y.~Feng, E.~W. Moore,
  J.~{VanderPlas}, D.~Laxalde, J.~Perktold, R.~Cimrman, I.~Henriksen, E.~A.
  Quintero, C.~R. Harris, A.~M. Archibald, A.~H. Ribeiro, F.~Pedregosa, P.~{van
  Mulbregt}, {SciPy 1.0 Contributors}, {{SciPy} 1.0: Fundamental Algorithms for
  Scientific Computing in Python}, Nature Methods 17 (2020) 261--272.

\bibitem{mitra2023guaranteed}
K.~Mitra, M.~Vohral{\'\i}k, {Guaranteed, locally efficient, and robust a
  posteriori estimates for nonlinear elliptic problems in iteration-dependent
  norms. An orthogonal decomposition result based on iterative linearization}
  (2023).

\bibitem{sun1998analysis}
H.~Sun, H.~Wong, K.~Ho, Analysis of infiltration in unsaturated ground, in:
  Proceedings of the annual seminar on slope engineering in Hong Kong, 1998,
  pp. 101--109.

\bibitem{lee2009simple}
L.~M. Lee, N.~Gofar, H.~Rahardjo, A simple model for preliminary evaluation of
  rainfall-induced slope instability, Engineering Geology 108~(3-4) (2009)
  272--285.

\bibitem{zhang2011stability}
L.~Zhang, J.~Zhang, L.~Zhang, W.~H. Tang, Stability analysis of
  rainfall-induced slope failure: a review, Proceedings of the Institution of
  Civil Engineers-Geotechnical Engineering 164~(5) (2011) 299--316.

\bibitem{bovolenta2020geotechnical}
R.~Bovolenta, D.~Bianchi, {Geotechnical Analysis and 3D Fem Modeling of Ville
  San Pietro (Italy)}, Geosciences 10~(11) (2020) 473.

\bibitem{report}
{Ministry for Primary Industries}, {ARPA Liguria - Rete di monitoraggio dei
  versanti (remover) regione Liguria, sito Ville San Pietro,},
  \url{https://srvcarto.regione.liguria.it/dtuff/img/Remover/Commenti_Siti/IM003_commento_tot.pdf}
  (2008-2022).

\bibitem{cressie1988spatial}
N.~Cressie, Spatial prediction and ordinary kriging, Mathematical geology 20
  (1988) 405--421.

\end{thebibliography}

\end{document}